\documentclass[11pt]{article}
\usepackage[table]{xcolor}
\usepackage{graphicx}
\usepackage{amsmath,amsfonts,amssymb}
\usepackage{graphicx}
\usepackage{multirow}
\usepackage{tikz,pgf}
\usepackage{url}
\usepackage{amsmath}
\usepackage{graphicx}
\usepackage{epsfig}
\usepackage{epstopdf}
\usepackage{color}
\usepackage{enumitem}
\def\tto{\;{\lower 1pt \hbox{$\rightarrow$}}\kern -10pt
\hbox{\raise 2pt \hbox{$\rightarrow$}}\;}
\def\Hat{\widehat}

\def\Bar{\overline}
\def\ra{\rangle}
\def\la{\langle}

\def\epsilon{\varepsilon}
\def\B{\Bbb B}
\def\h{\hfill\Box}
\def\R{\Bbb R}

\def\N{\Bbb N}
\def\ox{\bar{x}}

\def\oy{\bar{y}}
\def\oz{\bar{z}}

\def\ot{\bar{t}}

\def\co{\mbox{\rm co}\,}

\def\ri{\mbox{\rm ri}\,}

\def\gph{\mbox{\rm gph}\,}
\def\aff{\mbox{\rm aff}\,}
\def\epi{\mbox{\rm epi}\,}

\def\dom{\mbox{\rm dom}\,}
\def\aff{\mbox{\rm aff}\,}
\def\im{\mbox{\rm rge}\,}

\def\rge{\mbox{\rm rge}\,}

\def\h{\hfill\square}

\def\ph{\varphi}

\def\oR{\Bar{\R}}

\def\al{\alpha}

\def\ph{\varphi}

\def\oR{\Bar{\R}}

\def\al{\alpha}

\setlist[enumerate,1]{itemsep=0.0ex,parsep=0.5ex,label={\rm(\alph*)},leftmargin=*, align=left}
\newcounter{lk}

\usepackage[colorlinks=true]{hyperref}

\begin{document}
\begin{center}
{\sc\bf Near Convexity and Generalized Differentiation}\\[1ex]
{\sc Nguyen Mau   Nam}\footnote{Fariborz Maseeh Department of Mathematics and Statistics, Portland State University, Portland, OR
97207, USA (mnn3@pdx.edu). Research of this author was partly supported by the USA National Science Foundation under grant DMS-2136228.},
{\sc Nguyen Nang Thieu}\footnote{Institute of Mathematics, Vietnam Academy of Science and Technology,
Hanoi, Vietnam \& The State University of New York - SUNY, Korea (nnthieu@math.ac.vn).},
{\sc Nguyen Dong Yen}\footnote{Institute of Mathematics, Vietnam Academy of Science and Technology, 18 Hoang Quoc Viet,
Hanoi 10307 (ndyen@math.ac.vn).}\\
\end{center}
\small{\bf Abstract.} In this paper, we introduce the concept of nearly convex set-valued mappings and investigate fundamental properties of these mappings. Additionally, we establish a geometric approach for  generalized differentiation of nearly convex set-valued mappings and nearly convex functions. Our contributions expand the current knowledge of nearly convex sets and functions, while providing several new results pertaining to nearly convex set-valued mappings.\\[1ex]
{\bf Key words.}  Relative interior, nearly convex set, nearly convex function, nearly convex set-valued mapping,  subdifferential, normal cone, coderivative.\\[1ex]
\noindent {\bf AMS subject classifications.} 49J52, 49J53, 90C31

\newtheorem{Theorem}{Theorem}[section]
\newtheorem{Proposition}[Theorem]{Proposition}
\newtheorem{Remark}[Theorem]{Remark}
\newtheorem{Lemma}[Theorem]{Lemma}
\newtheorem{Corollary}[Theorem]{Corollary}
\newtheorem{Definition}[Theorem]{Definition}
\newtheorem{Example}[Theorem]{Example}
\renewcommand{\theequation}{\thesection.\arabic{equation}}
\normalsize

\section{Introduction}
\setcounter{equation}{0}

Over the past few decades, the classical notion of convexity and convex analysis have been extensively studied by prominent mathematicians and experts in applied fields. However, researchers have also made significant efforts to go beyond convexity by introducing and studying many generalized convexity notions for sets and functions. One such notion is the \textit{nearly convex set} introduced by Minty and Rockafellar, which is defined based on convexity by requiring that the set under consideration lies between a convex set and its closure in the Euclidean space $\R^n$; see~\cite{Minty1961, R1970}. It is worth noting that the concept of nearly convexity is sometimes referred to as \textit{almost convexity}. The study of maximally monotone operators reveals that the near convexity appears naturally, as their domains are nearly convex sets. In particular, the domain of the subdifferential mapping of a proper lower semicontinuous convex function is nearly convex, as discussed in~\cite{mmw2016} and the references therein.

Despite being introduced in the early age of convex analysis, the notion of near convexity had not been systematically studied in the literature until recently in~\cite{bmw2013,mmw2016}, which provide basic properties and examples of nearly convex sets,  ranges of maximally monotone operators, and ranges and fixed points of convex combinations of firmly nonexpansive mappings.  This development also creates opportunities for further study, from nearly convex sets to nearly convex functions and set-valued mappings. Progress in this direction can be seen in \cite{LM2019}, which builds on the work of Bo{\c{t}}, Kassay, and Wanka from \cite{bkw2008} on characterizations of nearly convex sets and functions, useful for investigating strong duality for nearly convex optimization problems. Further recent studies and applications of near convexity can be found in~\cite{nr1,nr2,nr3,nr4,nr5} and the references therein.

In this paper, we introduce the concept of nearly convex set-valued mappings by requiring their  graphs to  be nearly convex sets. We then examine how near convexity  of set-valued mappings and nonsmooth functions is preserved under different operations. Specifically, we present two distinct proofs that demonstrate how the sum of two nearly convex functions maintains near convexity, as long as the relative interiors of their effective domains intersect. Our proofs are straightforward and allowed us to identify an error in~\cite[Theorem~4.18]{LM2019}. Then we explore generalized differentiation for both nearly set-valued mappings and nonsmooth functions.

It should be noted that there are still challenging open questions related to nearly convex sets and functions. For example, to the best of our knowledge, the three questions raised by Moffat, Moursi, and Wang in ~\cite[p.~218]{mmw2016} on the domains and ranges of subdifferential mappings  have not been solved.

Our paper is structured as follows. In Section 2, we provide an overview of the fundamental definitions and properties of nearly convex sets, which will be referenced throughout the paper. In Section 3, we  establish basic properties of nearly convex functions and set-valued mappings. Section 4 is dedicated to exploring the preservation of near convexity under various operations on set-valued mappings and nonsmooth functions. In Section 5, we develop a geometric approach to generalized differentiation for nearly convex functions and set-valued mappings.

The paper utilizes standard notions and notations of convex analysis in the Euclidean space $\R^n$, which can be found in sources such as \cite{Borwein2000,HU2,bmn,r}. In the sequel, we use the notation $B(z,\rho)$ (resp., $\B (z,\rho)$) to represent the open (resp., closed) ball centered at $z\in\mathbb R^k$ with a radius of $\rho>0$. The affine hull of a subset $D\subset\mathbb R^k$ is abbreviated to $\aff D$, while $\overline{D}$ represents the closure of $D$, and $\mbox{\rm int}\,D$ denotes the interior of $D$.

\section{Preliminaries}
\setcounter{equation}{0}

This section provides a brief overview of fundamental definitions and significant properties of nearly convex sets that are utilized in this paper. For more comprehensive information, we recommend referring to~\cite{bmw2013,mmw2016}.

A subset $\Omega$ of $\R^n$ is said to be {\em nearly convex} if there exists a convex set $C$  such that
\begin{equation*}
C\subset \Omega\subset \Bar{C}.
\end{equation*}
Clearly, any convex set is nearly convex and any nearly convex subset of $\R$ is convex. Meanwhile, in $\R^n$ with $n\geq 2$ there are many nearly convex sets which are not convex.

Recall that the {\em relative interior} of an arbitrary set $\Omega$  in $\R^n$ is defined by
\begin{equation*}
\ri\Omega=\big\{a\in \Omega\; \big|\; \mbox{\rm there exists }\delta>0\; \mbox{\rm such that }B(a; \delta)\cap \aff\Omega\subset \Omega\big\}.
\end{equation*}
It follows from the definition that $a\in \ri\Omega$ if and only if $a\in \aff\Omega$ and there exists $\delta>0$ such that
$$B(a; \delta)\cap \aff\Omega\subset \Omega.$$
It is clear that if $\Omega_1\subset \R^n$ and $\Omega_2\subset \R^p$ are nearly convex, then $\Omega_1\times \Omega_2$ is also nearly convex.

Although the next proposition can be found in \cite{bmw2013}, we provide here a detailed proof for the convenience of the reader.
\begin{Proposition}\label{H} Let $\Omega$ be a nearly convex set with $C\subset \Omega\subset \Bar{C}$, where $C$ is a convex set in $\R^n$. Then $\aff\Omega=\aff C$, $\ri C=\ri\Omega$, and $\overline{C}=\overline{\Omega}$.
\end{Proposition}
{\bf Proof.} It follows from the definition that
\begin{equation*}
\aff C\subset \aff\Omega\subset \aff\overline{C}=\aff C,
\end{equation*}
which implies the equality $\aff\Omega=\aff C$. Fix any $a\in \ri C$. Then  $a\in C$ and there exists $\delta>0$ such that
\begin{equation*}
B(a; \delta)\cap \aff \Omega=B(a; \delta)\cap \aff C\subset C\subset \Omega.
\end{equation*}
It follows that $a\in \ri \Omega$, so $\ri C\subset \ri \Omega$.  The same argument shows that $\ri\Omega\subset \ri\overline{C}=\ri C$, where the last equality is valid by the convexity of $C$. Thus we obtain $\ri\Omega=\ri C$.

Since $C\subset \Omega\subset \overline{C}$, we have
\begin{equation*}
\overline{C}\subset\overline{\Omega}\subset\overline{C},
\end{equation*}
which implies $\overline{\Omega}=\overline{C}$ and completes the proof. $\h$

\begin{Theorem}\label{T1}  {\rm (See~\cite[Theorem~4.2 and Corollary~4.8]{mmw2016})} Suppose that $\Omega,\Omega_1,\dots, \Omega_m$ are nearly convex sets in $\R^n$, and $A\colon \R^n\to \R^p$ is a linear function. Then
\begin{enumerate}
\item $A(\Omega)$ is a nearly convex set in $\R^p$ and $\ri A(\Omega)=A(\ri\Omega)$.
\item If $\bigcap\limits_{i=1}^m\ri\Omega_i\neq\emptyset$, then $\bigcap\limits_{i=1}^m\Omega_i$ is nearly convex and
\begin{equation*}
\ri\left(\bigcap_{i=1}^m\Omega_i\right)=\bigcap_{i=1}^m\ri\Omega_i.
\end{equation*}
\end{enumerate}
\end{Theorem}

For two subsets $\Omega_1$ and $\Omega_2$ of $\R^n$,  if $\ri\Omega_1=\ri\Omega_2$ and $\overline{\Omega_1}=\overline{\Omega_2}$, then we say that they are {\em nearly equal} and write $\Omega_1\approx\Omega_2$.

The proposition below characterizes the near convexity using the near equality.

\begin{Proposition}\label{T2} {\rm (See \cite[Lemma~2.9]{bmw2013})} Let $\Omega$ be a subset of $\R^n$. Then the following properties are equivalent:
	\begin{enumerate}
		\item $\Omega$ is nearly convex.
		\item $\Omega$ is nearly equal to a convex set.
		\item $\Omega$ is nearly equal to a nearly convex set.
		\item $\Omega\approx \co\Omega$.
	\end{enumerate}
\end{Proposition}

The preservation of the near equality between two nearly convex sets via a linear mapping  is stated as follows.
\begin{Proposition}\label{T2a}  {\rm (See~\cite[Corollary~4.9]{mmw2016})}
	Let $\Omega_1$ and $\Omega_2$ be nearly convex subsets of $\R^n$ with $\Omega_1 \approx \Omega_2$ and $A\colon \R^n\to \R^p$ be a linear function. Then $A(\Omega_1) \approx A(\Omega_2)$.
\end{Proposition}

Given a function $f\colon \R^n\to \oR:=[-\infty, \infty]$, the {\em effective domain} and the  {\em epigraph} of $f$ are given  respectively by
\begin{align*}
&\dom f=\{x\in \R^n\; |\; f(x)<\infty\},\\
&\epi f=\{(x, \lambda)\in \R^n\times \R\; |\; f(x)\leq \lambda\}.
\end{align*}
We say that $f$ is {\em proper} if $\dom f\neq\emptyset$ and $-\infty <f(x)$ for all $x\in \R^n$. Throughout the paper we deal mostly with proper functions but occasionally encounter improper ones. The function $f$ is said to be {\em convex} if $\epi f$ is a convex set, and it is said to be {\em nearly convex} if $\epi f$ is a nearly convex set.

\begin{Example}{\rm {\bf (a)} Consider the function $f\colon \R\to \oR$ defined by
\begin{equation*}
f(x)=\begin{cases}
-\infty\; &\mbox{\rm if }x<0,\\
0\; &\mbox{\rm if }x\geq 0.
\end{cases}
\end{equation*}
Then $\dom f=\R$, $\epi f=\big((-\infty, 0)\times \R\big)\cup \big([0, \infty)\times [0, \infty)\big)$. We can see that $f$ is neither convex nor nearly convex.\\[1ex]
{\bf (b)} Consider the function $g\colon \R\to \oR$ defined by
\begin{equation*}
g(x)=\begin{cases}
-\infty\; &\mbox{\rm if }x<0,\\
1\; &\mbox{\rm if }x=0,\\
\infty&\mbox{\rm if }x>0.
\end{cases}
\end{equation*}
Then $\dom f=(-\infty, 0]$, $\epi f=\big((-\infty, 0)\times \R\big)\cup \big(\{0\}\times [ 1, \infty)\big)$. We can see that $f$ is an improper convex function.
}
\end{Example}

We continue this section with a representation of the affine hull  of the epigraph of an arbitrary function; see \cite[Exercise~2.6]{bmn}.

\begin{Proposition}\label{T3} Let $f\colon \R^n\to \oR$ be  a proper function. Then
\begin{equation}\label{aff epi}
\mbox{\rm aff}(\epi f)=\mbox{\rm aff}(\dom f)\times \R.
\end{equation}
\end{Proposition}
{\bf Proof.} Fix any $(x, \gamma)\in \mbox{\rm aff}(\epi f)$ and find $\lambda_i\in \R$ and $(x_i, \gamma_i)\in \epi f$ for $i=1, \ldots, m$  such that $\sum\limits_{i=1}^m\lambda_i=1$ and
\begin{equation*}
(x, \gamma)=\sum\limits_{i=1}^m \lambda_i(x_i, \gamma_i).
\end{equation*}
Since $(x_i, \gamma_i)\in \epi f$, we have $f(x_i)\leq \gamma_i<\infty$, so $x_i\in \dom f$ for $i=1, \ldots, m$. Then $x=\sum\limits_{i=1}^m \lambda_i x_i\in \mbox{\rm aff}(\dom f)$ and thus $(x, \gamma)\in \mbox{\rm aff}(\dom f)\times \R$. This justifies the inclusion~$\subset$ in~\eqref{aff epi}.

To verify the reverse inclusion in \eqref{aff epi}, take an arbitrary element $(x, \gamma)\in \mbox{\rm aff}(\dom f)\times \R$. Then $\gamma\in \R$ and there exist $\lambda_i\in \R$ for $i=1, \ldots, m$ with $\sum\limits_{i=1}^m\lambda_i=1$ such that $x=\sum\limits_{i=1}^m\lambda_i x_i$. Define $\alpha_i=f(x_i)\in \R$ for $i=1, \ldots, m$ and let $\alpha=\sum\limits_{i=1}^m \lambda_i \alpha_i\in \R$.  Clearly, $(x_i, \alpha_i)\in \epi f$ and $(x_i, \alpha_i+1)\in \epi f$. It follows that
\begin{align*}
&\sum_{i=1}^m \lambda_i(x_i, \alpha_i)=(x, \alpha)\in \mbox{\rm aff}(\epi f),\\
&\sum_{i=1}^m \lambda_i(x_i, \alpha_i+1)=(x, \alpha+1)\in \mbox{\rm aff}(\epi f).
\end{align*}
Considering the number $\lambda=\alpha-\gamma+1$, we have
\begin{equation*}
(x, \gamma)=\lambda (x, \alpha)+(1-\lambda)(x, \alpha+1)\in \mbox{\rm aff}(\epi f),
\end{equation*}
which justifies the reverse inclusion in~\eqref{aff epi} and also the proof of the  proposition. $\h$

For a set-valued mapping $F\colon \R^n\tto \R^p$, define the {\em domain},  the \textit{range}, and the {\em graph} of $F$ by
\begin{align*}
&\dom F=\{x\in \R^n\; |\; F(x)\neq\emptyset\},\quad  \rge F=\bigcup\limits_{x\in\R^n}F(x),\\
&\gph F=\{(x, y)\in \R^n\times \R^p\; |\; y\in F(x)\}.
\end{align*}
We say that $F$ is  {\em nearly convex} if $\gph F$ is a nearly convex set in $\R^n\times \R^p$.

Given a proper function $f\colon \R^n\to \oR$, define the {\em epigraphical mapping} $E_f\colon \R^n\tto \R$ by
\begin{equation}\label{Emapping}
   E_f(x)=[f(x), \infty)=\big\{\lambda\in \R\;\big|\; f(x)\leq \lambda\big\},\ \; x\in \R^n.
\end{equation}
It follows directly from the definition that $\dom E_f=\dom f$ and $\gph E_f=\epi f$. We also define the {\em epigraphical range} of $f$ by $\rge f=\rge E_f$.

\section{Nearly Convex Functions and Set-Valued Mappings}
\setcounter{equation}{0}

In this section, we study general properties of nearly convex functions and set-valued mappings. In particular, we are able to show that the relative interior of a nearly convex set-valued mapping can be represented in terms of the relative interior of its domains as well as those of its mapping values. This result generalizes a well-known theorem by Rockafellar on relative interiors of convex graphs to the case of nearly convex graphs;  see  \cite[Theorem~6.8]{r} and also \cite[Proposition~2.43]{rw} for another proof.

For two elements $a, b\in \R^n$, define
\begin{align*}
&[a, b]=\{(1-t)a+tb\; |\; 0\leq t\leq 1\},\\
&(a, b)=\{(1-t)ta+tb\; |\; 0<t<1\},\\
&[a, b)=\{(1-t)a+tb\; |\; 0\leq t< 1\}.
\end{align*}
Note that if $a=b$, then $[a,b]=(a,b)=[a, b)=\{a\}=\{b\}$.

The next simple result will be used in what follows.

\begin{Proposition}\label{Connection} Let $\Omega$ be a nearly convex set in $\R^n$. If $a\in \ri\Omega$ and $b\in \overline{\Omega}$, then $$[a, b)\subset \ri\Omega.$$
\end{Proposition}
{\bf Proof.} Take  any $a\in \ri\Omega$ and $b\in \overline{\Omega}$. Choose a convex set $C$ such that
$C\subset \Omega\subset \overline{C}$. By Proposition~\ref{H} we have $a\in \ri C$ and $b\in \overline{C}$. It follows that
$$[a, b)\subset \ri C=\ri\Omega,$$
which completes the proof. $\h$

Given two nonempty sets $\Omega_1$ and $\Omega_2$,  we say that $\Omega_1$ and $\Omega_2$ can be {\em properly separated} (by a hyperplane) if there exits $v\in \R^n$ such that the following two inequalities are satisfied:
\begin{eqnarray}\label{propsep}
\begin{array}{ll}
&\sup\big\{\la v, x\ra\;\big |\; x\in \Omega_1\big\}\leq \inf\big\{\la v, y\ra\;\big |\; y\in \Omega_2\big\},\\
&\inf\big\{\la v, x\ra\;\big |\; x\in \Omega_1\big\}< \sup\big\{\la v, y\ra\;\big |\; y\in \Omega_2\big\}.
\end{array}
\end{eqnarray}

Note that the first inequality means that $\la v, x\ra\leq \la v, y\ra$ whenever $x\in \Omega_1$ and $y\in \Omega_2$, while the second inequality means that there exist $\Hat{x}\in \Omega_1$ and $\Hat{y}\in \Omega_2$ such that $\la v, \Hat{x}\ra<\la v, \Hat{y}\ra$.

The theorem below provides necessary and sufficient conditions for proper separation of two nearly convex sets; see~\cite[Proposition~3.7]{LM2019}. Here we give a new proof for the result.

\begin{Theorem}\label{main sep} Let $\Omega_1$ and $\Omega_2$ be two nonempty nearly convex sets in $\R^n$. Then $\Omega_1$ and $\Omega_2$ can be properly separated if and only if $\ri\Omega_1\cap \ri\Omega_2=\emptyset.$
\end{Theorem}
{\bf Proof.} Let $C_1$ and $C_2$ be two convex sets  in $\R^n$ such that
\begin{equation}\label{new_label_inclusions}
C_1\subset \Omega_1\subset \Bar{C_1}\ \; \mbox{\rm and }\; C_2\subset \Omega_2\subset \Bar{C_2}.
\end{equation}
Suppose that $\Omega_1$ and $\Omega_2$ can be properly separated and find $v\in \R^n$ such that~\eqref{propsep} is satisfied.
The inclusions  in~\eqref{new_label_inclusions} and the continuity of the inner product yield
\begin{eqnarray*}
\begin{array}{ll}
\sup\big\{\la v, x\ra\;\big |\; x\in C_1\big\}&\leq \sup\big\{\la v, x\ra\;\big |\; x\in \Omega_1\big\}\leq \inf\big\{\la v, y\ra\;\big |\; y\in \Omega_2\big\}\\
&\leq \inf\big\{\la v, y\ra\;\big |\; y\in C_2\big\}.
\end{array}
\end{eqnarray*}
and
\begin{eqnarray*}
\begin{array}{ll}
\inf\big\{\la v, x\ra\;\big |\; x\in C_1\big\}&=\inf\big\{\la v, x\ra\;\big |\; x\in \Bar{C_1}\big\}\leq \inf\big\{\la v, x\ra\;\big |\; x\in \Omega_1\big\}\\
&<\sup\big\{\la v, y\ra\;\big |\; y\in \Omega_2\big\}\leq \sup\big\{\la v, y\ra\;\big |\; y\in \Bar{C_2}\big\}\\
&=\sup\big\{\la v, y\ra\;\big |\; y\in C_2\big\}.
\end{array}
\end{eqnarray*}
Since $C_1$ and $C_2$ are nonempty and convex, we can apply~\cite[Theorem~2.40]{bmn} and get that $\ri C_1\cap \ri C_2=\emptyset$.  By Proposition~\ref{H} we have
\begin{equation*}
\ri\Omega_1\cap \ri\Omega_2=\ri C_1\cap \ri C_2=\emptyset.
\end{equation*}
For the converse implication suppose that $\ri \Omega_1\cap \ri\Omega_2=\emptyset$ and get from Proposition~\ref{H} that $\ri C_1\cap \ri C_2=\emptyset$. Applying \cite[Theorem~2.40]{bmn} again gives us a vector $v\in \R^n$ such that
\begin{eqnarray}\label{new_label_separation}
\begin{array}{ll}
&\sup\big\{\la v, x\ra\;\big |\; x\in C_1\big\}\leq \inf\big\{\la v, y\ra\;\big |\; y\in C_2\big\},\\
&\inf\big\{\la v, x\ra\;\big |\; x\in C_1\big\}< \sup\big\{\la v, y\ra\;\big |\; y\in C_2\big\}.
\end{array}
\end{eqnarray}
Then  by~\eqref{new_label_inclusions} and~\eqref{new_label_separation} we have
\begin{eqnarray*}
\begin{array}{ll}
\sup\big\{\la v, x\ra\;\big |\; x\in \Omega_1\big\}\leq \sup\big\{\la v, x\ra\;\big |\; x\in \Bar{C_1}\big\}&=\sup\big\{\la v, x\ra\;\big |\; x\in C_1\big\}\\
&\leq \inf\big\{\la v, y\ra\;\big |\; y\in C_2\big\}\\
& = \inf\big\{\la v, y\ra\;\big |\; y\in\Bar{C_2}\big\}\\
&\leq \inf\big\{\la v, y\ra\; \big|\; y\in \Omega_2\big\}.
\end{array}
\end{eqnarray*}
The verification of the strict inequality $\inf\{\la v, x\ra\; |\; x\in \Omega_1\}<\sup\{\la v, y\ra\; |\; y\in \Omega_2\}$  by using~\eqref{new_label_separation} and ~\eqref{new_label_inclusions} is similar. Thus~\eqref{propsep} is satisfied for the vector $v$. This completes the proof. $\h$

\begin{Proposition}
	\label{prop:ri_ext}Let $\Omega$ be a nearly convex set in $\R^n$ with $y_0\in\Omega$. Then $y_0\in \ri\Omega$ if and only if for any $x\in \Omega$ there exists $z\in \Omega$ such that $y_0\in (x, z)$.
\end{Proposition}
{\bf Proof.} Suppose that $y_0\in \ri\Omega$ and take any $x\in \Omega$. By the definition  of relative interior, there exists $\delta>0$ such that
\begin{equation*}
B(y_0; \delta)\cap \aff\Omega\subset\Omega.
\end{equation*}
Choose $t>0$ sufficiently small such that $y_0+t(y_0-x)=(1+t)y_0+(-t)x\in B(y_0; \delta)$. Since $z=y_0+t(y_0-x)$ is an affine combination of $y_0$ and $x$, we see that $u\in B(y_0; \delta)\cap \aff\Omega\subset\Omega$. Then
\begin{equation*}
y_0=\frac{t}{1+t} x+\frac{1}{1+t}z\in (x, z).
\end{equation*}

To prove the converse implication, suppose on the contrary that for any $x\in \Omega$ there exists $z\in \Omega$ such that $y_0\in (x, z)$, but $y_0\notin \ri\Omega$. Choose a convex set $C$ such that
$ C\subset \Omega\subset \overline{C}.$ Then $\ri C=\ri \Omega$ and $\overline{C}=\overline{\Omega}$  by Proposition~\ref{H}; thus $y_0\notin \ri C$. Applying the separation theorem (see, e.g., \cite[Theorem~2.40]{bmn}) to the convex sets  $C$ and $\{y_0\}$, we can find $v\in \R^n$ such that
\begin{equation}\label{sep in}
\la v, x\ra\leq \la v, y_0\ra\ \; \mbox{\rm for all }\; x\in C
\end{equation}
and there exists $\hat{x}\in C$ such that $\la v,\hat{x}\ra<\la v, y_0\ra$. Since $\hat{x}\in \Omega$,  by our assumption there exists $\hat{z}\in \Omega$ such that  $y_0\in (\hat{x}, \hat{z})$. So, thanks to the inclusion $\Omega\subset \overline{C}$ and the convexity of $\overline{C}$, we have $y_0\in (\hat{x}, \hat{z})\subset [\hat{x}, \hat{z}]\subset \overline{C}$. Passing to a limit shows that  the inequality in~\eqref{sep in} holds for all $x\in \overline{C}$. Then
\begin{equation*}
\la v, x\ra\leq \la v, y_0\ra\ \; \mbox{\rm for all }\; x\in [\hat{x}, \hat{z}]
\end{equation*}
and $\la v, \hat{x}\ra <\la v, y_0\ra$.  This means that the convex sets $[\hat{x}, \hat{z}]$ and $\{y_0\}$ can be properly separated. So, by \cite[Theorem~2.40]{bmn} we obtain
\begin{equation*}
y_0\notin \mbox{\rm ri}([\hat x, \hat z])=(\hat x, \hat z),
\end{equation*}
which is a contradiction. $\h$

\medskip
Recall that a function is said to be nearly convex if its epigraph is nearly convex.

\begin{Proposition}\label{nc dom} If $F\colon \R^n\tto \R^p$ is a nearly convex set-valued mapping, then $\dom F$  and $\im F$ are a nearly convex sets. Consequently, if  $f\colon\R^n\to \oR$ is proper and nearly convex, then $\dom f$ and $\rge f$ are both nearly convex.
\end{Proposition}
{\bf Proof.} First observe that  $\dom F=\mathcal{P}(\gph F)$, where $\mathcal{P}$ is the linear mapping
\begin{equation}\label{Proj}
\mathcal{P}(x,y)=x, \ \; (x,y)\in \R^n\times \R^p.
\end{equation}
By Theorem \ref{T1}(a), the set $\dom F$ is nearly convex. Similarly, as  $\rge F=\mathcal{P}_1(\gph F)$, where $\mathcal{P}_1$ is the linear mapping
\begin{equation}\label{Proj_1}
	\mathcal{P}_1(x,y)=y, \ \; (x,y)\in \R^n\times \R^p.
\end{equation}
Now, suppose that $f$ is nearly convex. Then the epigraphical mapping $E_f$ defined in~\eqref{Emapping} is nearly convex. Since $\dom E_f=\dom f$ and $\rge E_f=\rge f$, the sets $\dom f$ and $\rge f$ are both nearly convex. $\h$

\begin{Remark}{\rm Given a proper function $f\colon \R^n\to \oR$, define
\begin{equation*}
\mbox{\rm im}\, f=\big\{f(x)\; \big |\; x\in \dom f\big\}.
\end{equation*}
We can find an example of a proper convex function for which $\mbox{\rm im}\,f$ is not nearly convex. Indeed, consider the function $f\colon \R\to \oR$ given by
\begin{equation*}
f(x)=\begin{cases}
1\; &\mbox{\rm if }x=0,\\
0\; &\mbox{\rm if }0<x\leq 1,\\
\infty&\mbox{\rm otherwise }.
\end{cases}
\end{equation*}
Then $f$ is convex, $\mbox{\rm im}\, f=\{0, 1\}$, and $\rge f=[0, \infty)$.
}\end{Remark}

\medskip
The theorem below allows us to represent the relative interior of the graph of a nearly convex set-valued mapping via the relative interiors of its domain and mapping values.

\begin{Theorem}\label{ri gph rep} Let $F\colon \R^n\tto \R^p$ be  a nearly convex set-valued mapping. Then we have
\begin{equation}\label{ri im}
\mbox{\rm ri}(\gph F)=\{(\ox, \oy)\in \R^n\times \R^p\; \big |\; \ox\in \mbox{\rm ri}(\dom F), \; \oy\in \ri F(\ox)\big\}.
\end{equation}
\end{Theorem}
{\bf Proof.} Since $\gph F$ is nearly convex by our assumption, using the projection mapping $\mathcal{P}$ defined in~\eqref{Proj} along with Theorem~\ref{T1}(a) gives us
\begin{equation}\label{ri proj}
\mathcal{P}(\mbox{\rm ri}(\gph F))=\ri\big(\mathcal{P}(\gph F)\big)=\mbox{\rm ri}(\dom F).
\end{equation}

To prove the inclusion $\subset$ in~\eqref{ri im}, fix any $(\ox, \oy)\in \mbox{\rm ri}(\gph F)$. By~\eqref{ri proj}
we have $\ox\in \mbox{\rm ri}(\dom F)$. By  the definition of relative interior, there exists $\delta>0$ such that
\begin{equation*}
\big[\B(\ox; \delta)\times \B(\oy; \delta)\big]\cap \mbox{\rm aff}(\gph F)\subset \gph F.
\end{equation*}
Then we have
\begin{equation}\label{ri gph2}
\big[\{\ox\}\times \B(\oy; \delta)\big]\cap \mbox{\rm aff}(\gph F)\subset \gph F,
\end{equation}
which implies that
\begin{equation}\label{ri gph3}
\B(\oy; \delta)\cap \aff F(\ox)\subset F(\ox).
\end{equation}
Indeed, taking any $y\in \B(\oy; \delta)\cap \aff F(\ox)$ gives us the representation
\begin{equation*}
y=\sum_{i=1}^m \lambda _i y_i,
\end{equation*}
where $y_i\in F(\ox)$ and $\lambda_i\in\R$ for $i=1, \ldots, m$ with $\sum\limits_{i=1}^m \lambda_i=1$. Then
\begin{equation*}
(\ox, y)=\sum_{i=1}^m \lambda_i(\ox, y_i)\in \mbox{\rm aff}(\gph F).
\end{equation*}
Thus, by~\eqref{ri gph2} we see that
\begin{equation*}
(\ox, y)\in \big[\{\ox\}\times \B(\oy; \delta)\big]\cap \mbox{\rm aff}(\gph F)\subset \gph F,
\end{equation*}
which implies that $y\in F(\ox)$. This justifies~\eqref{ri gph3}, so $\oy\in \ri F(\ox)$ by the definition of relative interior.

To prove the inclusion $\supset$ in~\eqref{ri im}, take any $\ox \in \mbox{\rm ri}(\dom F)$ and $\oy\in \ri F(\ox)$. Using~\eqref{ri proj}, we find $\hat{y}\in F(\ox)$ such that $(\ox, \hat{y})\in \mbox{\rm ri}(\gph F)$. Choose a convex set $C\subset \R^n\times \R^p$ such that
\begin{equation*}
C\subset \gph F\subset \Bar{C}.
\end{equation*}
Then, by Proposition~\ref{H} one has ${\rm ri}(\gph F)=\ri C$. We only need to consider the case where $\oy\neq \hat{y}$ because in the other case it holds that $(\ox,\oy)\in \mbox{\rm ri}(\gph F)$. Since $\oy\in \ri F(\ox)$, by \cite[Proposition~2.18]{bmn},  we can choose $\gamma>0$ such that $\oy+\gamma(\oy-\hat{y})\in F(\ox)$. Note that this does not require the convexity or near convexity of $F(\ox)$. Then $(\ox, \oy+\gamma(\oy-\hat{y}))\in \gph F\subset \Bar{C}$. By~\cite[Theorem~2.22]{bmn} and the convexity of $C$, we have the inclusion
\begin{equation*}
\big[ (\ox, \hat{y}), (\ox, \oy+\gamma(\oy-\hat{y}))\big)\subset \ri C,
\end{equation*}
where the set on the left-hand side is the half-open interval connecting $(\ox, \hat{y})\in \ri C$ with $(\ox, \oy+\gamma(\oy-\hat{y}))$. It follows that
\begin{equation}\label{interval}
	\big[ (\ox, \hat{y}), (\ox, \oy+\gamma(\oy-\hat{y}))\big)\subset {\rm ri}(\gph F).\end{equation}
Choosing $t=1/(1+\gamma)\in (0, 1)$, by~\eqref{interval} we obtain
\begin{equation*}
(\ox, \oy)=(1-t)(\ox, \hat{y})+ t\big(\ox, \oy+\gamma(\oy-\hat{y})\big)\in \big[ (\ox, \hat{y}), (\ox, \oy+\gamma(\oy-\hat{y}))\big)\subset \mbox{\rm ri}(\gph F),
\end{equation*}
which completes the proof of the theorem. $\h$

\medskip
The proposition below not only improves Proposition~4.4 from~\cite{LM2019} but also provides an alternative simple proof for the result.

 \begin{Proposition}\label{riepi} Let $f\colon \R^n\to \oR$ be a proper function. Then
 \begin{equation}\label{epi rep1}
 \mbox{\rm ri}(\epi f)\subset\big\{(x, \lambda)\in \R^n\times \R\;\big |\; x\in \mbox{\rm ri}(\dom f),\ \lambda>f(x)\big\}.
 \end{equation}
 The reverse inclusion of~\eqref{epi rep1} holds if we assume in addition that $f$ is nearly convex.
  \end{Proposition}
{\bf Proof.} To prove the first assertion, take any $(\ox, \overline{\lambda})\in \mbox{\rm ri}(\epi f)$. Then by the definition of relative interior $(\ox, \overline{\lambda})\in \epi f$, and there exists $\delta>0$ such that
\begin{equation*}
	\big (\B(\ox; \delta)\times [\overline{\lambda}-\delta, \overline{\lambda}+\delta]\big)\cap \mbox{\rm aff}(\epi f) \subset \epi f.
\end{equation*}
The representation of $\mbox{\rm aff}(\epi f)$ from Proposition~\ref{T3} gives us
\begin{equation*}
	\big (\B(\ox; \delta)\times [\overline{\lambda}-\delta, \overline{\lambda}+\delta]\big)\cap \big(\mbox{\rm aff}(\dom f)\times \R\big)=\big(\mbox{\rm aff}(\dom f)\cap \B(\ox; \delta)\big)\times [\overline{\lambda}-\delta, \overline{\lambda}+\delta]\subset \epi f.
\end{equation*}
It follows that
\begin{align*}
	&\big(\mbox{\rm aff}(\dom f)\cap \B(\ox; \delta)\big)\times\{\overline{\lambda}\}\subset \epi f,\\
	&\{\ox\}\times [\overline{\lambda}-\delta, \overline{\lambda}+\delta]\subset\epi f.
\end{align*}
The first inclusion gives us $\mbox{\rm aff}(\dom f)\cap \B(\ox; \delta)\subset \dom f$, and so $\ox\in \mbox{\rm ri}(\dom f)$. From the second inclusion we have $(\ox, \overline{\lambda}-\delta)\in \epi f$, so $f(\ox)\leq \overline{\lambda}-\delta<\overline{\lambda}$. This shows that $(\ox, \overline{\lambda})$ is contained in the set on the right-hand side of inclusion~\eqref{epi rep1}. Thus,~\eqref{epi rep1} is valid.

Now, to prove the second assertion, assume that $f$ is nearly convex. Consider the epigraphical mapping $F=E_f$ defined in~\eqref{Emapping}. Take $\ox\in \mbox{\rm ri}(\dom f)$ and $\overline{\lambda}>f(\ox)$. We see that $$\ox\in \mbox{\rm ri}(\dom f)=\mbox{\rm ri}(\dom  E_f)=\mbox{\rm ri}(\dom F).$$ Obviously, $\bar{\lambda}\in (f(\ox), \infty)={\rm ri}\big(E_f(\ox)\big)={\rm ri}\big(F(\ox)\big)$. Thus, $$(\bar{x}, \bar{\lambda})\in\mbox{\rm ri}(\gph F) =\mbox{\rm ri}(\gph E_f)=\mbox{\rm ri}(\epi f)$$ by the reverse inclusion in~\eqref{ri im}. $\h$

Given a function $f\colon \R^n\to \oR$, recall  that $f$ is continuous at $\ox\in \R^n$ if $\ox\in \mbox{\rm int}(\dom f)$ and for any $\epsilon>0$ there exists $\delta>0$ such that
\begin{equation}\label{C def}
    |f(x)-f(\ox)|<\epsilon\ \; \mbox{\rm whenever }\; x\in \B(\ox; \delta)\subset \dom f.
\end{equation}

The next proposition provides representations for the interior of the graph of a nearly convex set-valued mapping and also the interior of the epigraph of a nearly convex function.

\begin{Proposition}\label{intepi} Let $F\colon \R^n\tto \R^p$ be a set-valued mapping with $\mbox{\rm int}(\gph F)\neq\emptyset$. Then we have the inclusion
\begin{equation}\label{int gph}
\mbox{\rm int}(\gph F)\subset\big\{(\ox, \oy)\in \R^n\times \R^p\; \big|\; \ox\in \mbox{\rm int}(\dom F), \; \oy\in \mbox{\rm int}\, F(\ox)\big\}.
\end{equation}
The reverse inclusion in~\eqref{int gph} holds if we assume in addition that $F$ is nearly convex. Consequently, if a function $f\colon \R^n\to \oR$  is proper, nearly convex, and continuous at some point $\ox\in \R^n$, then
\begin{equation}\label{int epi}
    \mbox{\rm int}(\epi f)=\big\{(x, \lambda)\in \R^n\times \R\; \big |\; x\in \mbox{\rm int}(\dom f), \; f(x)<\lambda\big\}.
\end{equation}
 \end{Proposition}
{\bf Proof}. Fix any $(\ox, \oy)\in \mbox{\rm int}(\gph F)$ and find $\delta>0$ such that
\begin{equation*}
\B(\ox; \delta)\times \B(\oy; \delta)\subset \gph F.
\end{equation*}
Then it holds that $\B(\ox; \delta)\times \{\oy\}\subset \gph F$ and $\{\ox\}\times \B(\oy; \delta)\subset \gph F$, which imply that  $\ox\in \mbox{\rm int}(\dom F)$ and $\oy\in \mbox{\rm int}\, F(\ox)$. Thus we have inclusion $\subset$ in~\eqref{int gph}. The reverse inclusion follows from Theorem~\ref{ri gph rep} under the assumption that $F$ is nearly convex.

Now, let $f\colon \R^n\to \oR$  be proper, nearly convex, and continuous at some $\ox\in \R^n$. It follows from the definition of continuity that $\ox\in \mbox{\rm int}(\dom f)$, so $\mbox{\rm ri}(\dom f)=\mbox{\rm int}(\dom f).$
Fix $\epsilon>0$ and choose $\delta>0$ such that \eqref{C def} is satisfied. Then
\begin{equation*}
\B(\ox; \delta)\times (f(\ox)+\epsilon, \infty)\subset \epi f,
\end{equation*}
so $\mbox{\rm int}(\epi f)\neq \emptyset$. Thus, representation \eqref{int epi} follows directly from \eqref{int gph} using again the epigraphical mapping $F=E_f$. $\h$

\medskip
A fundamental property of the values of a nearly convex set-valued mapping is given in the following theorem.

\begin{Theorem}
		\label{cor:nearly_convex} Let $F\colon \R^n\tto \R^p$ be a nearly convex set-valued mapping. If $\ox\in \mbox{\rm ri}(\dom F)$, then $F(\ox)$ is nearly convex. In particular, $\ri F(\ox)$ is nonempty.
	\end{Theorem}
{\bf Proof.} Suppose that $F$ is nearly convex and $\ox\in \mbox{\rm ri}(\dom F)$. Then $\gph F$ is nearly convex. Let $C=\mbox{\rm co}(\gph F)$ and define the set-valued mapping
	\begin{equation*}
		G(x)=\big\{y\in \R^p\; \big|\; (x, y)\in C\big\}, \ \; x\in \R^n.
\end{equation*}
Then $\gph G = C$. Since $F$ is nearly convex and $C=\mbox{\rm co}(\gph F)$, one has $C\approx \gph F$ by Proposition \ref{T2}. Let us consider the mapping $\mathcal{P}$ defined in \eqref{Proj}. By Proposition~\ref{T2a} we get $\mathcal{P}(\gph G) \approx \mathcal{P}(\gph F)$ or equivalently, $\dom G\approx\dom F$. Thus, $\mbox{\rm ri}(\dom G) =\mbox{\rm ri}(\dom F)$. Take any $\oy\in \ri G(\ox)$. Then by Theorem~\ref{ri gph rep} one has $(\ox,\oy)\in\mbox{\rm ri}(\gph G)= \mbox{\rm ri}(\gph F)$, which implies that $\oy\in \ri F(\ox)\subset F(\ox)$. Thus, we have
\begin{equation}
	\label{eq:includ_1}
	\ri G(\ox) \subset F(\ox)
\end{equation} In addition, we can deduce from the convexity of $G$ that $G(\ox)$ is convex. Since $\gph F\subset\gph G$, one has
\begin{equation}
	\label{eq:includ_2}
	F(\ox) \subset G(\ox)\subset \overline{G(\ox)}.
\end{equation}
By~\cite[Proposition~2.33]{bmn}, one obtains $\overline{\ri G(\ox)} = \overline{G(\ox)}$. Then the inclusions~\eqref{eq:includ_1} and~\eqref{eq:includ_2} yield that $F(\ox)$ is nearly convex. Since $F(\ox)$ is nearly convex and nonempty, we can deduce from Proposition~\ref{H} that $\ri F(\ox)\neq\emptyset$. $\h$

\medskip
If $F\colon \R^n\tto \R^p$ is a nearly convex set-valued mapping, then $F(x)$ is not necessarily nearly convex for all $x\in\dom F$.

\begin{Example}\label{Example1} {\rm
		Let $F\colon\R\tto \R$ be given by
		\begin{equation*}
				F(x) = \begin{cases}
			[0,2]\quad&\text{if}\; x\in [0,1),\\
			[0,1)\cup (1,2] & \text{if}\; x = 1,\\
			\emptyset &\text{otherwise.}
			\end{cases}
		\end{equation*}
	It is clear to see that $F$ is nearly convex, but $F(1)$ is not nearly convex.}
\end{Example}

A natural question arises: \textit{Whether under the assumption of Theorem~\ref{cor:nearly_convex} can one assert that $F(\ox)$ is convex for every $\ox\in \mbox{\rm ri}(\dom F)$?} As shown by the example below, the answer is negative.

\begin{Example}\label{Example2} {\rm
	Consider the constant set-valued mapping $F\colon\R\tto \R^2$ with $$F(x)=\big([0,1]\times [0,1]\big)\setminus \{(1/2,1)\},\ \; x\in\R.$$ Then $F$ is nearly convex, $\dom F=\R$, and $F(x)$ is not convex for all $x\in\R$.}
\end{Example}

The next theorem shows that the notions of convexity and near convexity are the same for functions defined on the real line.

\begin{Proposition}
		 If a proper function $f\colon \R\to \oR$ is nearly convex, then it is convex.
\end{Proposition}
{\bf Proof.} Suppose that $f\colon \R\to \oR$ is proper and nearly convex. First, observe that for any $x\in \mbox{\rm ri}(\dom f)$, $u\in \dom f$, and $0<t<1$ we have
\begin{equation}\label{JS}
f(tx+(1-t)u)\leq tf(x)+(1-t)f(u).
\end{equation}
Indeed, taking $\epsilon>0$ yields $(x, f(x)+\epsilon)\in \mbox{\rm ri}(\epi f)$ by Proposition~\ref{riepi}. So, applying Proposition~\ref{Connection} to the nearly convex set $\Omega=\epi f$ and noting that $(u, f(u))\in \epi f$, we have $t(x, f(x)+\epsilon)+(1-t)(u, f(u))\in \mbox{\rm ri}(\epi f)$. Using Proposition~\ref{riepi} again gives us
\begin{equation*}
f(tx+(1-t)u)<tf(x)+(1-t)f(u)+t\epsilon,
\end{equation*}
which implies \eqref{JS} by letting $\epsilon\to 0^+$. Now, we continue the proof with the observation that  $\dom f\subset \R$ is nearly convex by Proposition \ref{nc dom}. It is easy to see that a nearly convex subset of $\R$ is an interval or a singleton, so $\dom f$ is an interval  in $\R$ or a singleton. The function $f$ is clearly convex if its domain is a singleton. It suffices to consider the case where $\dom f = [\alpha,\beta]$ with $\alpha, \beta\in \R$ and  $\alpha <\beta$ because the conclusion is obvious for other cases. By the observation above, we only need to prove that for any $t\in (0,1)$ one has
			\begin{equation*}
				f(t\al +(1-t)\beta) \leq t f(\al)+(1-t) f(\beta).
			\end{equation*}
			Take any  $0<t<1$ and let $\hat{x} = t\al +(1-t)\beta$. Since $\hat{x}\in \mbox{\rm ri}(\dom f)$, one has $\frac{\al+\hat{x}}{2}\in \mbox{\rm ri}(\dom f)$ and $\frac{\beta+\hat{x}}{2}\in\mbox{\rm ri}(\dom f)$.
We can deduce from~\eqref{JS} that
			\begin{equation*}
			f\left(\frac{\al+\hat{x}}{2}\right)\leq \frac{f(\al)+f(\hat{x})}{2}
\quad\text{and}\quad
			f\left(\frac{\beta+\hat{x}}{2}\right)\leq \frac{f(\beta)+f(\hat{x})}{2}.
		\end{equation*}
So, using~\eqref{JS} again gives us
	\begin{align*}
	f(\hat{x})=	f\left(t\frac{\al+\hat{x}}{2}+ (1-t)\frac{\beta+\hat{x}}{2}\right)&\leq 	t f\left(\frac{\al+\hat{x}}{2}\right) + (1-t)f\left(\frac{\beta+\hat{x}}{2}\right)\\
		&\leq t \frac{f(\al)+f(\hat{x})}{2} + (1-t)\frac{f(\beta)+f(\hat{x})}{2}\\
		&=\frac{t f(\al)+(1-t)f(\beta)}{2} +\frac{f(\hat{x})}{2}.
	\end{align*}
This yields
\begin{equation*}
f(\hat{x}) \leq \frac{t f(\al)+(1-t)f(\beta)}{2}.
\end{equation*}
 Therefore, $f$ is convex as desired. $\h$

\medskip
A nearly convex function $f\colon \R^n\to \oR$, $n\geq 2$, needs not to be convex.

\begin{Example}{\rm Consider the function $f\colon \R^2\to \oR$ defined by
\begin{equation*}
f(x,y)=\begin{cases}
0\; &\mbox{\rm if } (x, y)\in \big([-1, 1]\times [-1, 1]\big)\setminus \{(1, 0)\},\\
\infty\; &\mbox{\rm otherwise}.
\end{cases}
\end{equation*}
Then $f$ is nearly convex but not convex.
}
\end{Example}

Given a set $\Omega\subset \R^n\times \R$, define $f_\Omega\colon \R^n\to [-\infty, \infty]$ by
\begin{equation*}
f_\Omega(x)=\inf\{t\in \R\; |\; (x, t)\in \Omega\},\; x\in \R^n.
\end{equation*}
Let $f\colon \R^n\to \oR$ be a function. Define $\co f\colon \R^n\to [-\infty, \infty]$ by $\co f=f_{\co(\epi f)}$, i.e.,
\begin{equation*}
(\co f)(x)=\inf\{t\in \R\; |\; (x, t)\in \mbox{\rm co}(\epi f)\}, \; x\in \R^n.
\end{equation*}
It follows from the definition that $\co f$ is the largest convex function majorized by $f$, i.e.,
\begin{align*}
&\co f \; \mbox{\rm is convex},\\
&\co f \leq f,\\
&\mbox{\rm if }g\colon \R^n\to [-\infty, \infty] \; \mbox{\rm is a convex function such that }g\leq f, \; \mbox{\rm the }g\leq\co f.
\end{align*}

Recall that two subsets $\Omega_1$ and $\Omega_2$ of $\R^n$ are nearly equal and we write $\Omega_1\approx\Omega_2$ if $\ri\Omega_1=\ri\Omega_2$ and $\overline{\Omega}_1=\overline{\Omega_2}$. We say that two functions $f_1, f_2\colon \R^n\to \oR$ are nearly equal and write $f_1\approx f_2$ if their epigraphs are nearly equal sets.
\begin{Lemma}\label{coepi} Let $f\colon \R^n\to \oR$ be a proper function. Then we have the inclusions
\begin{equation*}
\mbox{\rm co}(\epi f)\subset \mbox{\rm epi}(\co f)\subset \overline{\mbox{\rm co}(\epi f)}.
\end{equation*}
\end{Lemma}
{\bf Proof.} Since $\epi f\subset \mbox{\rm epi}(\co f)$, where the latter set is a convex due to the convexity of $\co f$, we see that $\mbox{\rm co}(\epi f)\subset \mbox{\rm epi}(\co f)$.

Next, observe that $\mbox{\rm co}(\epi f)$ is an {\em epigraphical set} in the sense that if $(x_0, t_0)\in \mbox{\rm co}(\epi f)$, then $(x_0, t_0+\alpha)\in \mbox{\rm co}(\epi f)$ for every $\alpha\geq 0$. Indeed, if $(x_0, t_0)\in \mbox{\rm co}(\epi f)$, then we can write
\begin{equation*}
(x_0, t_0)=\sum_{i=1}^m \lambda_i(x_i, t_i),
\end{equation*}
where $(x_i, t_i)\in \epi f$, $\lambda_i\geq 0$ for all $i=1, \ldots, m$, and $\sum_{i=1}^m \lambda_i=1$. It follows that
\begin{equation*}
(x_0, t_0+\alpha)=\sum_{i=1}^m \lambda_i(x_i, t_i+\alpha).
\end{equation*}
Since $(x_i, t_i)\in \epi f$ and $\alpha\geq 0$, one has
\begin{equation*}
(x_i, t_i+\alpha)\in \epi f\; \mbox{\rm for all }i=1, \ldots, m.
\end{equation*}
It follows that $(x_0, t_0+\alpha)\in \mbox{\rm co}(\epi f)$ as desired.

Now, take $(\ox, \overline{t})\in \mbox{\rm epi}(\co f)$ and suppose on the contrary that $(\ox, \overline{t})\notin \overline{\mbox{\rm co}(\epi f)}$. Choose $\rho>0$ and $\epsilon>0$ such that
\begin{equation*}
\big[B(\ox; \rho)\times [\ot-\epsilon, \ot+\epsilon]\big]\cap \overline{\mbox{\rm co}(\epi f)}=\emptyset.
\end{equation*}
In particular, $\big[(\{\ox\}\times [\ot-\epsilon, \ot+\epsilon]\big]\cap \overline{\mbox{\rm co}(\epi f)}=\emptyset$.

We have $\ot\geq (\co f)(\ox)$, and by the definition there exists a sequence $\{t_k\}\subset \R$ such that $(\ox, t_k)\in \mbox{\rm co}(\epi f)$ and $t_k\to (\co f)(\ox)$ as $k\to \infty$. Consider the case where $(\co f)(\ox)\in \R$. In this case, for some $k_0\in \N$ we have
\begin{equation*}
(\co f)(\ox)\leq t_k<(\co f)(\ox)+\epsilon\; \mbox{\rm when }k\geq k_0.
\end{equation*}
{\bf Case 1}: $(\co f)(\ox)\geq \ot-\epsilon$.  In this case $\ot-\epsilon\leq t_k\leq \ot+\epsilon$ for all $k\geq k_0$. Thus
\begin{equation*}
\big[B(\ox; \rho)\times [\ot-\epsilon, \ot+\epsilon]\big]\cap \overline{\mbox{\rm co}(\epi f)}\neq\emptyset
\end{equation*}
because $(\ox, t_k)\in \mbox{\rm co}(\epi f)\subset \overline{\mbox{\rm co}(\epi f)}$ for all $k\geq k_0$. This yields a contradiction. \\[1ex]
{\bf Case 2}: $(\co f)(\ox)< \ot-\epsilon$. In this case since $t_k\to (\co f)(\ox)$, there exists $\hat{k}\in \N$ such that
\begin{equation*}
t_k<\ot-\epsilon<\ot\; \mbox{\rm for all }k\geq \hat{k}.
\end{equation*}
Fix $k>\hat{k}$. Since $(\ox, t_k)\in \mbox{\rm co}(\epi f)$ and $t_k<\ot$, by the claim above we see that $(\ox, \ot)\in \mbox{\rm co}(\epi f)\subset \overline{\mbox{\rm co}(\epi f)}$. This again yields a contradiction.

For the case where $(\co f)(\ox)=-\infty$, we have $t_k\to-\infty$ as $k\to \infty$. Then there exists  $\hat{k}\in\N$ such that
\begin{equation*}
t_k<\ot\; \mbox{\rm for all }k\geq \hat{k}.
\end{equation*}
For any $k\geq \hat{k}$, since $(\ox, t_k)\in \mbox{\rm co}(\epi f)$ and $t_k<\ot$, we obtain $(\ox, \ot)\in \mbox{\rm co}(\epi f)\subset \overline{\mbox{\rm co}(\epi f)}$ by the epigraphical property of $\mbox{\rm co}(\epi f)$. This is a contradiction to the choice of $(\ox, \ot)$.
$\h$

The following proposition can be found in \cite[Theorem~4.12]{LM2019} for which we provide here a self-contained proof for the convenience of the reader.

\begin{Proposition} \label{prop:funct_nearly_equal}Let $f\colon \R^n\to \oR$ be a proper function. Then $f$ is nearly convex if and only if $f\approx \co f$.
\end{Proposition}
{\bf Proof.} Suppose that $f\approx \co f$. Since $\co f$ is a convex function, its epigraph $\mbox{\rm epi}(\co f)$ is a convex set in $\R^n\times \R$. By the definition of near equality we have $\epi f\approx \mbox{\rm epi}(\co f)$, which means that $\epi f$ is nearly equal to a convex set. Thus, by Proposition~\ref{T2} $\epi f$ is nearly convex, and so $f$ is nearly convex.

For the converse implication, suppose that $f$ is a nearly convex function. Then $\epi f$ is a nearly convex set in $\R^n\times\R$. By the definition, there exists a convex set $C$ in $\R^{n+1}$ such that
\begin{equation*}
C\subset \epi f\subset \overline{C}.
\end{equation*}
Since $\overline{C}$ is closed and convex, we see that
\begin{equation*}
C\subset \epi f\subset \mbox{\rm co}(\epi f)\subset \overline{\mbox{\rm co}(\epi f)}\subset \overline{C}.
\end{equation*}
By Lemma \ref{coepi} we have $\mbox{\rm co}(\epi f)\subset \mbox{\rm epi}(\co f)\subset \overline{\mbox{\rm co}(\epi f)}$, which yields
\begin{equation*}
C\subset \epi f\subset \mbox{\rm co}(\epi f)\subset \mbox{\rm epi}(\co f ))\subset \overline{\mbox{\rm co}(\epi f)}\subset \overline{C}.
\end{equation*}
This implies by Proposition~\ref{H} that $\epi f\approx \mbox{\rm epi}(\co f)$ and therefore $f\approx \co f$. $\h$

\section{Preservation of Near Convexity under Basic Operations}
\setcounter{equation}{0}

This section focuses on studying the preservation of near convexity in set-valued mappings and nonsmooth functions under  basic operations on them. Our research presents significant advancements beyond the recent findings in \cite{LM2019}. Notably, we offer two distinct proofs demonstrating the preservation of near convexity under the summation of two nearly convex functions. We also provide a counterexample to illustrate that this result does not apply unless a relative interior qualification condition is assumed.

\begin{Proposition}\label{nearly convexity of sum} Let $f\colon \R^n\to \oR$ be a proper nearly convex function. Define the function $\ph\colon \R^n\times \R\to \oR$ by
\begin{equation*}
\ph(x, \alpha)=f(x)+\alpha,\ \; (x, \alpha)\in \R^n\times \R.
\end{equation*}
Then $\ph$ is also nearly convex.
\end{Proposition}
{\bf Proof.} Define the sets
\begin{align*}
&\Theta_1=\epi \ph=\{(x, \alpha, \lambda)\in \R^n\times \R\times \R\; |\; f(x)+\alpha\leq \lambda\},\\
&\Theta_2=\{(x, \alpha, \lambda)\in \R^n\times\R\times \R\; |\; (\co f)(x)+\alpha\leq\lambda\}.
\end{align*}
Then $\Theta_2$ is a convex set. To verify that $\ph$ is nearly convex, it suffices to show that $\Theta_1=\epi \ph$ is nearly convex. To accomplish this goal, by Proposition~\ref{T2} we will show that $\Theta_1\approx\Theta_2$.\\[1ex]
{\bf Step 1.} $\overline{\Theta}_1=\overline{\Theta}_2$: Since $\co f\leq f$, we see that $\Theta_1\subset\Theta_2$ and thus $\overline{\Theta}_1\subset\overline{\Theta}_2$. To verify the reverse inclusion, let us show that $\Theta_2\subset\overline{\Theta}_1$. Fix any $(x, \alpha, \lambda)\in \Theta_2$. Then $(\co f)(x)+\alpha \leq \lambda$, so $(\co f)(x)\leq \lambda-\alpha$. It follows that
\begin{equation*}
(x, \lambda-\alpha)\in \mbox{\rm epi}(\co f)\subset \overline{\mbox{\rm epi}(\co f)}=\overline{\epi f}
\end{equation*}
due to the fact that $f\approx \co f$ from Proposition~\ref{prop:funct_nearly_equal}. Thus there exists a sequence $\{(x_k, \gamma_k)\}\subset \epi f$ such that
\begin{equation*}
(x_k, \gamma_k)\to (x, \lambda-\alpha)\; \mbox{\rm as }k\to \infty.
\end{equation*}
Then $f(x_k)\leq \gamma_k$, so $f(x_k)+\alpha\leq \gamma_k+\alpha$ and thus $(x_k, \alpha, \gamma_k+\alpha)\in\Theta_1$. Since $(x_k, \alpha, \gamma_k+\alpha)$ converges to $(x, \alpha, \lambda)$, this yields $(x, \alpha, \lambda)\in \overline{\Theta}_1$ and completes the proof of Step 1.\\[1ex]
{\bf Step 2.} $\ri\Theta_1=\ri\Theta_2$: Since $\overline{\Theta}_1=\overline{\Theta}_2$, we see that $\aff\Theta_1=\aff\Theta_2$, which yields $\ri\Theta_1\subset\ri\Theta_2$ due to the fact that $\Theta_1\subset\Theta_2$. It follows from Proposition~\ref{T3} that
\begin{equation*}
\aff \Theta_1=\mbox{\rm aff}(\dom f)\times \R\times\R.
\end{equation*}
By Proposition~\ref{prop:funct_nearly_equal} one has $\mbox {\rm epi}(f) \approx \mbox{\rm epi}(\co f)$. Using the projection mapping, it follows from Proposition~\ref{T2a} that $\dom f \approx \mbox{\rm dom}(\co f)$. Thus, $\overline{\dom f} = \overline{ \mbox{\rm dom}(\co f)}$ and hence,
\begin{equation*}
	\mbox{\rm aff}(\dom f) =  \mbox{\rm aff}(\overline{\dom f}) =\mbox{\rm aff}(\overline{ \mbox{\rm dom}(\co f)}) = \mbox{\rm aff}( \mbox{\rm dom}(\co f)).
\end{equation*}
Then we have
\begin{equation*}
	\aff \Theta_2=\mbox{\rm aff}(\dom f)\times \R\times\R =\aff\Theta_1.
\end{equation*}  By Proposition \ref{riepi} we have
\begin{align*}
\ri \Theta_2&=\{(x, \alpha, \lambda)\; |\; x\in \mbox{\rm ri}(\mbox{\rm dom}(\co (f))), (\co f)(x)+\alpha<\lambda\}\\
&=\{(x, \alpha, \lambda)\; |\; x\in \mbox{\rm ri}(\mbox{\rm dom}(\co f)), (\co f)(x)<\lambda-\alpha\}\\
&=\{(x, \alpha, \lambda)\; |\; (x, \lambda-\alpha)\in \mbox{\rm ri} (\mbox{\rm epi}(\co f))\}\\
&=\{(x, \alpha, \lambda)\; |\; (x, \lambda-\alpha)\in \mbox{\rm ri} (\epi f)\}\\
&=\{(x, \alpha, \lambda)\; |\; x\in \mbox{\rm ri}(\dom f), f(x)<\lambda-\alpha\}.
\end{align*}
Take any $(x_0, \alpha_0, \lambda_0)\in \ri \Theta_2$ and get $(x_0, \lambda_0-\alpha_0)\in \mbox{\rm ri}(\epi f)$. Thus by the definition of the relative interior and Proposition~\ref{T3} we find $\gamma>0$ and $\epsilon >0$ such that
\begin{equation*}
\big[B(x_0; \gamma)\times (\lambda_0-\alpha_0-\epsilon, \lambda_0-\alpha_0+\epsilon)\big]\cap \big[\aff(\dom f)\times \R\big]\subset \epi f.
\end{equation*}
Choose $\delta>0$ such that if $(\alpha, \lambda)\in (\alpha_0-\delta, \alpha_0+\delta)\times (\lambda_0-\delta, \lambda_0+\delta)$, then $\lambda-\alpha\in (\lambda_0-\alpha_0-\epsilon, \lambda_0-\alpha_0+\epsilon)$ (use the continuity of $g(\alpha, \lambda)=\lambda-\alpha$). Now we see that
\begin{equation*}
\big[B(x_0; \gamma)\times (\alpha_0-\delta, \alpha_0+\delta)\times (\lambda_0-\delta, \lambda_0+\delta)\big]\cap \big[\mbox{\rm aff}(\dom f)\times \R\times\R\big]\subset \Theta_1.
\end{equation*}
Indeed, if $(x, \alpha, \lambda)$ is in the set on the left-hand side of this inclusion. Then
\begin{equation*}
(x, \lambda-\alpha)\in \big[B(x_0; \gamma)\times (\lambda_0-\alpha_0-\epsilon, \lambda_0-\alpha_0+\epsilon)\big]\cap \big[\mbox{\rm aff}(\dom f)\times \R\big]\subset \epi f.
\end{equation*}
Thus $f(x)\leq \lambda-\alpha$, and so $(x, \alpha, \lambda)$ is in the set on the right-hand side of the inclusion. This shows that $(x_0, \alpha_0, \lambda_0)\in \ri\Theta_1$. Therefore, $\ri\Theta_2\subset \ri \Theta_1$. $\h$.

%

\begin{Theorem}
		\label{thm:sum_rule_setvalued}
Let $F_1, F_2\colon \R^n\tto \R^p$ be nearly convex set-valued mappings. Then $F_1+F_2$ is also nearly convex under the qualification condition
\begin{equation}\label{QCF}
	\mbox{\rm ri}(\dom F_1)\cap \mbox{\rm ri}(\dom F_2)\neq\emptyset.
\end{equation}
	\end{Theorem}
{\bf Proof.} Define two sets
\begin{eqnarray}\label{O12}
	\begin{array}{ll}
		&\Omega_1=\{(x, y_1, y_2)\in \R^n\times \R^p\times \R^p\; |\; y_1\in F_1(x)\} = (\gph F) \times \R^p,\\
		&\Omega_2=\{(x, y_1, y_2)\in \R^n\times \R^p\times \R^p\; |\; y_2\in F_2(x)\}.
	\end{array}
\end{eqnarray}
Since $F_1$ and $F_2$ are nearly convex, $\Omega_1$ and $\Omega_2$ are nearly convex. It follows from  Theorem \ref{ri gph rep} that
\begin{eqnarray*}
	\begin{array}{ll}
		&\ri\Omega_1=\{(x, y_1, y_2)\in \R^n\times \R^p\times \R^p\; |\; x\in \mbox{\rm ri}(\dom F_1), \; y_1\in \ri F_1(x)\},\\
		&\ri \Omega_2=\{(x, y_1, y_2)\in \R^n\times \R^p\times \R^p\; |\; x\in \mbox{\rm ri}(\dom F_1), \; y_2\in \ri F_2(x)\in \mbox\}.
	\end{array}
\end{eqnarray*}
We can choose  $\hat{x}\in \mbox{\rm ri}(\dom F_1)\cap \mbox{\rm ri}(\dom F_2)$ due to  \eqref{QCF}. By Theorem \ref{cor:nearly_convex} we can choose $\hat{y}_1\in \ri F_1(\ox)$ and $\hat{y}_2\in \ri F_2(\ox)$.  Thus, $(\hat{x}, \hat{y}_1, \hat{y}_2)\in \ri\Omega_1\cap \ri\Omega_2$.  Since $\Omega_1$ and $\Omega_2$ are nearly convex, it follows from Theorem~\ref{T1}(b) that $\Omega_1\cap \Omega_2$ is nearly convex.

Define the linear mapping $A\colon \R^n\times \R^p\times \R^p\to \R^n\times\R^p$ by
	\begin{equation*}
		A(x,y_1,y_2) = (x,y_1+y_2), \ \;  (x,y_1,y_2)\in \R^n\times \R^p\times \R^p.
	\end{equation*}
Obviously, $\mbox{\rm gph}(F_1+F_2) = A(\Omega_1\cap\Omega_2)$. It follows from Theorem~\ref{T1}(a)  that $\mbox{\rm gph}(F_1+F_2)$ is nearly convex, so $F_1+F_2$ is a nearly convex set-valued mapping.
$\h$

As a consequence, we obtain the corollary below. We not only provide a correct statement for \cite[Theorem 4.18]{LM2019} but also give a simple proof for this result.

\begin{Corollary}\label{nc sum f}  Let $f_1, f_2\colon \R^n\to \oR$ be proper nearly convex functions. Suppose that
\begin{equation}\label{QCf}
\mbox{\rm ri}(\dom f_1)\cap \mbox{\rm ri}(\dom f_2)\neq \emptyset.
\end{equation}
Then $f_1+f_2$ is nearly convex.
\end{Corollary}
{\bf Proof.} Consider the epigraphical  mappings $F_i=E_{f_i}$, $i=1, 2$, given in \eqref{Emapping}. Then $\dom F_i=\dom f_i$, $\gph F_i=\epi f_i$, and $(F_1+F_2)(x)=[f_1(x)+f_2(x), \infty)$ for all $x\in \R^n$. By \eqref{QCf}, we can choose $\hat{x}\in \mbox{\rm ri}(\dom f_1)\cap \mbox{\rm ri}(\dom f_2)$. Then choose a real number $\hat{\lambda}$ such that
\begin{equation*}
\max\{f_1(\hat{x}, f_2(\hat{x})\}<\hat{\lambda}.
\end{equation*}
It follows from Proposition \ref{riepi} that
\begin{equation*}
(\hat{x}, \hat{\lambda})\in \mbox{ri}(\epi f_1)\cap \mbox{\rm ri}(\epi f_2)=\mbox{\rm ri}(\gph F_1)\cap \mbox{\rm ri}(\gph F_2).
\end{equation*}
By Theorem \ref{thm:sum_rule_setvalued}, the mapping $F_1+F_2$ is nearly convex. Thus $\mbox{\rm gph}(F_1+F_2)=\mbox{\rm epi}(f_1+f_2)$ is nearly convex. By the definition, $f_1+f_2$ is nearly convex. $\h$

\begin{Remark} {\rm We can prove Corollary \ref{nc sum f} by an alternative way as follows. Define the sets
\begin{eqnarray*}
\begin{array}{ll}
&\Omega_1=\{(x, \lambda, \alpha)\in \R^n\times \R\times \R\; |\; f_1(x)\leq \alpha\},\\
&\Omega_2=\{(x, \lambda, \alpha)\in \R^n\times \R\times \R\; |\; \lambda\geq f_2(x)+\alpha\}.
\end{array}
\end{eqnarray*}
Let us show that $\mbox{\rm epi}(f_1+f_2)=\mathcal{P}(\Omega_1\cap \Omega_2)$, where $\mathcal{P}$ is the projection mapping onto $\R^n\times \R$ (removing the last component). The set $\Omega_2$ is nearly convex by Proposition \ref{nearly convexity of sum}.
Take any $(x, \lambda, \alpha)\in \Omega_1\cap \Omega_2$. Then $f_1(x)\leq \alpha\leq \lambda-f_2(x)$ and thus $\lambda\geq f_1(x)+f_2(x)$. It follows that $(x, \lambda)\in \mbox{\rm epi}(f_1+f_2)$. This implies that
$\mathcal{P}(\Omega_1\cap \Omega_2)\subset \mbox{\rm epi}(f_1+f_2)$.

Now fix any $(x, \lambda)\in \mbox{\rm epi}(f_1+f_2)$ and get $f_1(x)+f_2(x)\leq \lambda$. This implies $f_1(x)\leq \lambda-f_2(x)$. Choose a real number $\alpha$ such that $f_1(x)\leq \alpha\leq \lambda-f_2(x)$. Then $(x, \lambda, \alpha)\in \Omega_1\cap \Omega_2$ and so
\begin{equation*}
(x,\lambda)=\mathcal{P}(x, \lambda, \alpha)\in \mathcal{P}(\Omega_1\cap \Omega_1),
\end{equation*}
which justifies the inclusion $\mbox{\rm epi}(f_1+f_2)\subset\mathcal{P}(\Omega_1\cap \Omega_2)$.

It follows from Proposition~ \ref{riepi} and Proposition~\ref{nearly convexity of sum}  that
\begin{eqnarray*}
\begin{array}{ll}
&\ri \Omega_1=\{(x, \lambda, \alpha)\in \R^n\times \R\times \R\; |\; x\in \mbox{\rm ri}(\dom f_1), f_1(x)<\alpha\},\\
&\ri\Omega_2=\{(x, \lambda, \alpha)\in \R^n\times \R\times \R\; |\; x\in  \mbox{\rm ri}(\dom f_2), \lambda> f_2(x)+\alpha\}.
\end{array}
\end{eqnarray*}
(We may think about $\Omega_2$ as the epigraph of the function $g(x, \alpha)=f_2(x)+\alpha$).

Choose $x_0\in \mbox{\rm ri}(\dom f_1)\cap \mbox{\rm ri}(\dom f_2)$. Then pick $\alpha_0>f_1(x_0)$ and pick $\lambda_0>f_2(x_0)+\alpha_0$. It is obvious that $(x_0, \lambda_0, \alpha_0)\in \ri\Omega_1\cap \ri\Omega_2$. Thus $\ri\Omega_1\cap \ri\Omega_2\neq\emptyset$, so $\Omega_1\cap \Omega_2$ is nearly convex. Therefore, $ \mbox{\rm epi}(f_1+f_2)$ is nearly convex. }

\end{Remark}

The example below shows that the result in \cite[Theorem~4.18]{LM2019} is not correct without assuming that the relative interiors of the domains of the functions involved intersect each other.

\begin{Example}{\rm Let $\Omega_1=([-1, 0]\times [-1, 1])\setminus \{(0,0)\}\subset \R^2$, and let $\Omega_2=([0, 1]\times [-1, 1])\setminus \{(0, 0)\}\subset \R^2$. Define $f_i(x)=\delta(x; \Omega_i)$ for $x\in \R^2$ and $i=1, 2$. Then $f_1+f_2=\delta(\cdot; \Omega_1\cap \Omega_2)$. It is clear that both $f_1$ and $f_2$ are nearly convex but their sum is not nearly convex.
}\end{Example}

The next theorem shows that the near convexity is preserved under compositions of set-valued mappings.

\begin{Theorem}
	\label{thm:chain_rule_setvalued}
	Let $F\colon \R^n\tto \R^p$ and $G\colon \R^p\tto \R^q$ be nearly convex set-valued mappings. Suppose that
\begin{equation}\label{QCC}
\mbox{\rm ri}(\rge F)\cap \mbox{\rm ri}(\dom G)\neq\emptyset.
\end{equation}
Then $G\circ F\colon \R^n\tto \R^q$ is nearly convex.
\end{Theorem}
{\bf Proof.} Define the sets
\begin{eqnarray}\label{O12c}
	\begin{array}{ll}
		&\Omega_1= (\gph F) \times \R^q \;\subset \R^n\times\R^p\times\R^q,\\
		&\Omega_2=\R^n\times (\gph G)\;\subset \R^n\times\R^p\times\R^q.
	\end{array}
\end{eqnarray}
Then $\Omega_1$ and $\Omega_2$ are nearly convex. By Theorem~\ref{ri gph rep} we have
\begin{eqnarray*}
	\begin{array}{ll}
		&\ri\Omega_1=\{(x, y,z)\in \R^n\times \R^p\times \R^q\; |\; x\in \mbox{\rm ri}(\dom F), y\in \ri F(x)\},\\
		&\ri \Omega_2=\{(x, y, z)\in \R^n\times \R^p\times \R^q\; |\; y\in \mbox{\rm ri}(\dom G), z\in \ri G(y)\}.
	\end{array}
\end{eqnarray*}
Choose $\hat{y}\in \mbox{\rm ri}(\rge F) \cap \mbox{\rm ri}(\dom G)$. Since $\hat{y}\in  \mbox{\rm ri}(\dom G)$, by Theorem Theorem \ref{cor:nearly_convex} the set $ \ri G(\hat{y})$ is nonempty, so we can choose $\hat{z}\in \ri G(\hat{y})$. Since $\hat{y}\in \mbox{\rm ri}(\rge F)  = \mbox{\rm ri}(\dom F^{-1})$ and $F^{-1}$ is also nearly convex, we can choose $\hat{x}\in \ri F^{-1}(\hat{y}) $. It can be easily seen that $\hat{x}\in \mbox{\rm ri}(\dom F)$ by considering the linear mapping $T\colon \R^p\times\R^n\to \R^n\times \R^p$, where $T(y,x)=(x,y)$ for  $(y,x)\in \R^p\times\R^n$. Then $(\hat{x},\hat{y},\hat{z}) \in \ri\Omega_1\cap\ri\Omega_2$. By Theorem~\ref{T1}(a) the intersection $\Omega_1\cap\Omega_2$ is nearly convex.  Define $A\colon \R^n\times\R^p\times\R^q\to \R^n\times\R^q$ by
\begin{equation*}
	A(x,y,z) = (x,z),\ \; (x,y,z)\in \R^n\times\R^p\times \R^q.
\end{equation*}
Obviously, $A(\Omega_1\cap\Omega_2) =\gph (G\circ F)$. Thus, $\gph(G\circ F)$ is nearly convex.
$\h$

Theorem \ref{thm:chain_rule_setvalued} allows us to derive a new result on the near convexity of the composition of a nearly convex function and an affine mapping.

\begin{Corollary} Let $B\colon \R^n\to \R^p$ be the affine mapping defined by
\begin{equation*}
B(x)=Ax+b\; \mbox{\rm for }x\in \R^n,
\end{equation*}
where $A\in \R^{p\times n}$ and $b\in \R^p$, and let $f\colon \R^p\to \oR$ be a nearly convex function. Then $f\circ B\colon \R^n\to \oR$ is nearly convex under the qualification condition
\begin{equation}\label{QCCf}
B(\R^n)\cap \ri(\dom f)\neq\emptyset.
\end{equation}
\end{Corollary}
{\bf Proof.} Let $F(x)=\{B(x)\}$ for $x\in \R^n$, and let $G(x)=[f(x), \infty)$. Then $\gph F=\gph B$, $\gph G=\epi f$, and
\begin{equation*}
(G\circ F)(x)=[(f\circ B)(x), \infty)\; \mbox{\rm for all }\; x\in \R^n.
\end{equation*}
Using \eqref{QCCf}, we can choose $\hat{x}\in \R^n$ such that $\hat{y}=B(\hat{x})\in \mbox{\rm ri}(\dom f)$. Then fix $\hat{\lambda}>f(\hat{y})$. Since $\gph B$ is an affine set,  $(\hat{x}, \hat{y})\in \mbox{\rm ri}(\gph F)=\gph B$. Using Proposition~\ref{riepi} gives us $(\hat{y}, \hat{\lambda})\in \mbox{\rm ri}(\epi f)=\mbox{\rm ri}(\gph G)$. By Theorem \ref{thm:chain_rule_setvalued}, the composition $G\circ F$ is nearly convex. Thus, $\mbox{\rm epi}(f\circ B)=\mbox{\rm gph}(G\circ F)$ is nearly convex, so  $f\circ B$ is nearly convex. $\h$

 Given a set-valued mapping $G\colon \R^n\tto \R^p$ with $\Omega\subset\R^n$ and $\Theta\subset\R^p$, recall that
\begin{align*}
&G(\Omega)=\bigcup_{x\in \Omega}G(x),\\
&G^{-1}(\Theta)=\big\{x\in \R^n\; \big|\; G(x)\cap \Theta\neq\emptyset\big\}.
\end{align*}
The corollary below shows that the near convexity of sets is preserved under nearly convex set-valued mappings under direct and inverse images.
	\begin{Corollary}\label{cor:ncpreserve}
		Let $G\colon \R^n\tto\R^p$ be a nearly convex set-valued mapping. Then we have the following assertions:
\begin{enumerate}
 \item If $\Omega\subset\R^n$  is a nearly convex  set, then $G(\Omega)$ is also a nearly convex set provided that
		\begin{equation*}
			\mbox{\rm ri}(\dom G)\cap \ri \Omega\neq \emptyset.
		\end{equation*}
\item If $\Theta\subset \R^p$ is a nearly convex set, the $G^{-1}(\Theta)$ is also a nearly convex set provided that
\begin{equation}\label{QCIN}
	\mbox{\rm ri}(\rge G)\cap \ri \Theta \neq \emptyset.
\end{equation}
\end{enumerate}
	\end{Corollary}
{\bf Proof.} Let us define a function $F\colon \R\tto \R^n$ by
\begin{equation*}
	F(x)=\begin{cases}
		\Theta\quad&\text{if}\; x =0\\
		\emptyset &\text{otherwise}.
	\end{cases}
\end{equation*}
Then $\gph F = \{0\}\times \Theta$, and hence $F$ is nearly convex. Note that $\mbox{\rm rge}(G\circ F)= \Theta$. Thus, $F\circ G$ is nearly convex. We see that
\begin{equation*}
	(G\circ F)(x)=\begin{cases}
		F(\Omega)\quad&\text{if}\; x =0\\
		\emptyset &\text{otherwise}.
	\end{cases}
\end{equation*}
Therefore, $\gph (F\circ G)=\{0\}\times F(\Theta)$. Define the linear mapping $A\colon \R\times \R^p\to \R^p$ by $A(x,y)=y$ for $(x, y)\in \R\times \R^p$. Since $F(\Omega) = A\big(\gph (G\circ F)\big)$, it follows from  Theorem~\ref{T1}(a) that $F(\Theta)$ is nearly convex.

We can prove the assertion (b) using assertion (a). Indeed, define the mapping $H\colon \R^p\tto \R^n$ by
\begin{equation*}
H(y)=G^{-1}(y)=\big\{x\in \R^n\; \big| \; y\in G(x)\big\}, \; y\in \R^p.
\end{equation*}
Then observe $H$ is nearly convex, $G^{-1}(\Theta)=H(\Theta)$ and $\dom H=\rge G$. Then it is  straightforward to complete the proof. $\h$

Let $F_i\colon \R^n\tto \R^p$ for $i=1, \ldots, m$ be set-valued mappings. Define
\begin{equation*}
	\big(\bigcap_{i=1}^mF_i\big)(x) = F_1(x)\cap \ldots \cap F_m(x),\ \; x\in \R^n.
\end{equation*}
In the next theorem we discuss the near convexity of the intersection mapping.
\begin{Theorem}\label{NCI}
	Let $F_i\colon \R^n\tto \R^p$ for $i=1, \ldots, m$ be nearly convex set-valued mappings. Suppose
	\begin{equation}\label{QCM}
		\bigcap_{i=1}^m \mbox{\rm ri}(\gph F_i) \neq \emptyset.
	\end{equation}
Then $\bigcap_{i=1}^mF_i$ is nearly convex.
\end{Theorem}
{\bf Proof.} It is obvious that
\begin{equation*}
	\mbox{\rm gph}(\bigcap_{i=1}^mF_i) = \bigcap_{i=1}^m \gph F_i.
\end{equation*}
By Theorem \ref{T1}(b), the set $\bigcap_{i=1}^m \gph F_i$ is nearly convex under the qualification condition \ref{QCM}. Therefore, $\bigcap_{i=1}^mF_i$ is nearly convex.
$\h$

To conclude this section, we discuss the near convexity of the maximum function.
Given functions $f_i\colon \R^n \to \oR$ for $i=1, \ldots, m$, define
\begin{equation}\label{MF}
f(x)=\max\{f_i(x)\; |\; i=1, \ldots, m\}, \; x\in \R^n.
\end{equation}

\begin{Corollary}\label{NCMf} Let $f_i\colon \R^n\to \oR$ for $i=1, \ldots, m$ be nearly convex functions. Suppose that
\begin{equation*}
\bigcap_{i=1}^m \mbox{\rm ri}(\dom f_i)\neq\emptyset.
\end{equation*}
Then the maximum function $f$ defined in \eqref{MF} is nearly convex.
\end{Corollary}
{\bf Proof.} Consider the epigraphical mappins $F_i(x)=E_{f_i}$ given in \eqref{Emapping} for $i=1, \ldots, m$, and let $F=\bigcap_{i=1}^mF_i$. Then we see that
\begin{equation*}
\epi f=\bigcap_{i=1}^m \epi f_i=\bigcap_{i=1}^m \gph F_i=\gph F.
\end{equation*}
Choosing $x_0\in \bigcap_{i=1}^m \mbox{\rm ri}(\dom f_i)$ and let $\lambda=f(x_0)+1$, it follows from Proposition~\ref{riepi} that
\begin{equation*}
(x_0, \lambda)\in \bigcap_{i=1}^m \mbox{\rm ri}(\epi f_i)=\bigcap_{i=1}^m \mbox{\rm ri}(\gph F_i).
\end{equation*}
By Theorem \ref{NCI}, the set-valued mapping $F$ is nearly convex. Thus, $\epi f=\gph F$ is nearly convex and therefore $f$ is nearly convex.
$\h$

\section{Nearly Convex Generalized Differentiation}
\setcounter{equation}{0}

In this section, we explore the topic of generalized differentiation for nearly convex set-valued mappings and nearly convex functions using a geometric approach that has proven successful in convex analysis. Previous work on this approach can be found in \cite{bmn,bmn2022}. We present new calculus rules for the coderivatives of sums, compositions, and maxima of nearly convex set-valued mappings, as well as related refinements for subdifferentials of nearly convex functions.

Given a nearly convex set $\Omega$ in $\R^n$ with $\ox\in \Omega$, define the normal cone to $\Omega$ at $\ox$ by
\begin{equation*}
N(\ox; \Omega)=\big\{v\in \R^n\; |\; \la v, x-\ox\ra\leq 0\ \, \mbox{\rm for all }\, x\in \Omega\big\}.
\end{equation*}

\begin{Theorem} \label{NCI} Let $\Omega_1$ and $\Omega_2$ be nearly convex sets such that
\begin{equation*}
\ri\Omega_1\cap \ri\Omega_2\neq \emptyset.
\end{equation*}
Then $\Omega_1\cap \Omega_2$ is nearly convex and
\begin{equation}\label{ncir}
N(\ox; \Omega_1\cap \Omega_2)=N(\ox; \Omega_1)+N(\ox; \Omega_2)\ \, \mbox{\rm for all }\, \ox\in \Omega_1\cap \Omega_2.
\end{equation}
\end{Theorem}
{\bf Proof.} The proof is straightforward based on \cite[Theorem~2.56]{bmn} with the use of Theorem \ref{main sep} and Proposition~\ref{riepi}. We provide the detailed proof here for the convenience of the reader.

Fix any $v\in N(\ox; \Omega_1\cap \Omega_2)$ and get
\begin{equation*}
\la v, x-\ox\ra\leq 0\; \mbox{\rm for all }x\in \Omega_1\cap \Omega_2.
\end{equation*}
Define the sets
\begin{align*}
&\Theta_1=\big\{(x, \lambda)\in \R^n\times \R\; \big |\; x\in \Omega_1,\;\lambda\leq \la v, x-\ox\ra\big\},\\
&\Theta_2=\Omega_2\times [0, \infty).
\end{align*}
We first claim that $\Theta_1$ is nearly convex and
\begin{equation}\label{rit}
\ri\Theta_1=\big\{(x, \lambda)\in \R^n\times \R\; \big|\; x\in \ri\Omega_1, \; \lambda<\la v, x-\ox\ra\big\}.
\end{equation}
Indeed, define the  set
\begin{equation*}
A=\big\{(x, \lambda)\in \R^n\times \R\; \big|\; \lambda\leq \la v, x-\ox\ra\big\}.
\end{equation*}
Then $A$ is a nonempty convex set and (prove it!)
\begin{equation*}
\ri A=\big\{(x, \lambda)\in \R^n\times \R\; \big|\; \lambda< \la v, x-\ox\ra\big\}.
\end{equation*}
Let $B=\Omega_1\times \R$ and see that $B$ is nearly convex with $\ri B=(\ri \Omega_1)\times \R$. Observe that $\Theta_1=A\cap B$. Choosing $x_0\in \ri\Omega_1$ and choosing $\lambda_0<\la v, x_0-\ox\ra$ give us $(x_0, \lambda_0)\in \ri A\cap \ri B$. By Theorem \ref{T1}(b), the set $\Theta_1=A\cap B$ is nearly convex and $\ri\Theta_1=\ri A\cap \ri B$, so \eqref{rit} is satisfied. Obviously, $\Theta_2$ is nearly convex with $\ri\Theta_2=(\ri \Omega_1)\cap (0, \infty)$. Thus, we can easily check that $\ri\Theta_1\cap \ri\Theta_2=\emptyset$ with proof by contradiction. By the proper separation from Theorem \ref{main sep}, there exists $(w, \gamma)\in \R^n\times \R$ such that
\begin{equation}\label{a}
\la w, x\ra +\gamma \lambda \leq \la w, y\ra+\gamma \beta\; \mbox{\rm whenever }(x, \lambda)\in \Theta_1, \; (y, \beta)\in \Theta_2.
\end{equation}
In addition, there exist $(x_0, \lambda_0)\in \Theta_1$ and $(y_0, \beta_0)\in \Theta_2$ such that
\begin{equation}\label{b}
\la w, x_0\ra +\gamma \lambda_0<\la w, y_0\ra+\gamma \beta_0.
\end{equation}
Using \eqref{a} with $x=\ox$, $\lambda=0$, $y=\ox$, $\beta=1$ gives $\gamma\geq 0$. If $\gamma=0$, then we can use \eqref{a} and \eqref{b} along with the definition to see that $\Omega_1$ and $\Omega_2$ can be properly separated, so $\ri\Omega_1\cap \ri\Omega_2=\emptyset$ by Theorem \ref{main sep}, a contradiction. Thus, $\gamma>0$.

Next, using \eqref{a} with $x=\ox$, $\lambda=0$, $y\in \Theta_2$, and $\beta=0$ gives us
\begin{equation*}
\la w, \ox\ra\leq \la w, y\ra\; \mbox{\rm for all }y\in \Omega_2,
\end{equation*}
This implies that $-w\in N(\ox; \Omega_2)$. Then we can use \eqref{a} with $x\in \Omega_1$, $\lambda =\la v, x-\ox\ra$, $y=\ox$, and $\beta=0$ to get
\begin{equation*}
\la w, x\ra+\gamma \la v, x-\ox\ra\leq \la w, \ox\ra\; \mbox{\rm for all }x\in \Omega_1.
\end{equation*}
Dividing both sides of this inequality by $\gamma$ and rearranging the terms, we have
\begin{equation*}
\la v+\dfrac{w}{\gamma}, x-\ox\ra \leq 0\; \mbox{\rm for all }x\in \Omega_1.
\end{equation*}
It follows that $v+\frac{w}{\gamma}\in N(\ox; \Omega_1)$, so
\begin{equation*}
v\in -\dfrac{w}{\gamma}+N(\ox; \Omega_1)\subset N(\ox; \Omega_2)+N(\ox; \Omega_2).
\end{equation*}
This justifies the inclusion $\subset$ in \eqref{ncir}, while the reverse inclusion can be proved by the definition. $\h$

 Let $F\colon \R^n\tto \R^p$ be a nearly convex set-valued mapping and let $(\ox,\oy)\in \gph F$. The {\em coderivative} of $F$ at $(\ox,\oy)$ is the set-valued mapping $D^*F(\ox,\oy)\colon \R^p\tto \R^n$ with the values
	\begin{equation*}\label{cod}
		D^*F(\ox,\oy)(v)=\big\{u\in \R^n\;\big|\;(u,-v)\in N\big((\ox,\oy);\gph F\big)\big\},\ \, v\in \R^p.
	\end{equation*}

Given  $(\ox,\oy)\in\mbox{\rm gph}(F_1+F_2)$, define the set
\begin{equation*}\label{S}
	S(\ox,\oy)=\big\{(\oy_1,\oy_2)\in \R^p\times \R^p\;\big|\;\oy=\oy_1+\oy_2,\;\oy_i\in F_i(\ox),\;i=1,2\big\}.
\end{equation*}
The theorem below provides a coderivative sum rule for nearly convex set-valued mappings.

\begin{Theorem}\label{CSR} Let $F_1, F_2\colon \R^n\tto \R^p$ be nearly convex set-valued mappings. Suppose that the qualification condition \eqref{QCF} is satisfied.  Then the equality
	\begin{equation*}
		D^*(F_1+F_2)(\ox, \oy)(v)=D^*F_1(\ox, \oy_1)(v)+D^*F_2(\ox, \oy_2)(v)
	\end{equation*}
holds for every $v\in \R^p$ and $(\oy_1, \oy_2)\in S(\ox, \oy)$, where $S$ is defined in~\eqref{S}.
\end{Theorem}
{\bf Proof.}  By Theorem \ref{thm:sum_rule_setvalued}, the set-valued mapping $F_1+F_2$ is nearly convex.  Fix any $(\oy_1, \oy_2)\in S(\ox, \oy)$ and $v\in \R^p$. Fix any \begin{equation*}\label{simple_inclusion} u\in D^*(F_1+F_2)(\ox,\oy)(v).\end{equation*}
Then we have the inclusion $(u,-v)\in N((\ox,\oy);\mbox{\rm gph}(F_1+F_2))$. Consider the sets $\Omega_i$, $i=1, 2$, defined in \eqref{O12}. By the definition we have
\begin{equation*}\label{inclusion_sum_rule}
	(u,-v,-v)\in N((\ox,\oy_1,\oy_2);\Omega_1\cap\Omega_2).
\end{equation*}
The proof of Theorem~\ref{thm:sum_rule_setvalued} tells us that $\ri\Omega_1\cap \ri\Omega_2\neq\emptyset$.
Then we can employ Theorem~\ref{NCI} and get
\begin{equation*}
	(u,-v,-v)\in N((\ox,\oy_1,\oy_2);\Omega_1\cap\Omega_2)=N((\ox,\oy_1,\oy_2);\Omega_1)+N((\ox,\oy_1,\oy_2);\Omega_2).
\end{equation*}
Therefore, the rest of the proof follows from that of \cite[Theorem~3.37]{bmn}.  $\h$

Let $f\colon \R^n\to \oR$ be a proper nearly convex function. We define the {\em subdifferential} of $f$ at $\ox\in \dom f$ by
\begin{equation*}
\partial f(\ox)=\big \{v\in \R^n\; \big|\; \la v, x-\ox\ra\leq f(x)-f(\ox)\ \, \mbox{\rm for all }\, x\in \R^n \big\}.
\end{equation*}

The proposition below allows us to represent the subdifferential of a nearly convex function via the coderivative of the epigraphical mapping \eqref{Emapping}.

\begin{Proposition}\label{CS} If $f\colon \R^n \to \oR$ is a proper nearly convex function, then
\begin{equation*}
D^*E_f(\ox, f(\ox))(1)=\partial f(\ox),
\end{equation*}
where $E_f$ is defined in \eqref{Emapping}.
\end{Proposition}
{\bf Proof.} By the definition,
\begin{equation*}
D^*E_f(\ox, f(\ox))(1)=\big\{v\in \R^n\; \big |\; (v, -1)\in N((\ox, f(\ox)); \epi f)\big\}.
\end{equation*}
Thus, taking any  $v\in D^*E_f(\ox, f(\ox))(1)$ gives us
\begin{equation}\label{c1}
\la v, x-\ox\ra-(\lambda-f(\ox))\leq 0\; \mbox{\rm whenver }(x, \lambda)\in \epi f.
\end{equation}
Using this inequality with $x\in \dom f$ and $\lambda=f(x)$ gives
\begin{equation}\label{d1}
\la v, x-\ox\ra\leq f(x)-f(\ox)\; \mbox{\rm for all }x\in \dom f,
\end{equation}
which implies that $v\in \partial f(\ox)$ since $f(x)=\infty$ if $x\notin \dom f$.

Now, suppose that $v\in \partial f(\ox)$ and get \eqref{d1}, which obviously implies \eqref{c1} since $f(x)\leq \lambda$ whenever $(x, \lambda)\in \epi f$. Thus, $v\in D^*E_f(\ox, f(\ox))(1)$, which completes the proof. $\h$

\begin{Proposition}\label{sub pro} Let $f\colon \R^n\to \oR$ be a proper nearly convex function. The following assertions hold:
\begin{enumerate}
    \item If $\ox\in \dom f$ and $(v, -\alpha)\in N((\ox, f(\ox)); \epi f)$, then $\alpha\geq 0$.
    \item If $\ox\in \mbox{\rm ri}(\dom f)$, then $\partial f(\ox)\neq\emptyset$. In particular, if $f$ is continuous at $\ox$, then $\partial f(\ox)\neq\emptyset$.
    \item If $f$ is continuous at $\ox$, then $N((\ox, \bar{\lambda}); \epi f)=\{(0, 0)\}$ whenever $f(\ox)<\bar{\lambda}$.
    \item $(v, 0)\in N((\ox, f(\ox)); \epi f)$ if and only if $v\in N(\ox; \dom f)$.
    \item If $\alpha>0$, then $(v, -\alpha)\in N((\ox, f(\ox)); \epi f)$ if and only if $v\in \alpha \partial f(\ox)$.
    \end{enumerate}
\end{Proposition}
{\bf Proof.} (a) It follows from the definition that
\begin{equation*}
\la v, x-\ox\ra-\alpha (\lambda -f(\ox))\leq 0\; \mbox{\rm whenever }f(x)\leq\lambda.
\end{equation*}
Using this inequality with $x=\ox$ and $\lambda=f(\ox)+1$ gives us the conclusion. \\[1ex]
(b) It follows from Proposition~\ref{riepi} that $(\ox, f(\ox))\notin \mbox{\rm ri}(\epi f)$. By the proper separation from Theorem \ref{main sep}, there exist $v\in \R^n$ and $\gamma\in \R$ such that
\begin{equation}\label{e}
\la v, x\ra-\gamma \lambda\leq \la v, \ox\ra-\gamma f(\ox)\; \mbox{\rm whenever }f(x)\leq \lambda.
\end{equation}
In addition, there exist $(x_0, \lambda_0)\in \epi f$ such that
\begin{equation*}
\la v, x_0\ra-\gamma \lambda_0< \la v, \ox\ra-\gamma f(\ox).
\end{equation*}
First, we see that $\gamma\geq 0$ by using \eqref{e} with $x=\ox$ and $\lambda=f(\ox)+1$. If $\gamma=0$, then we see that $\la v, x\ra\leq \la v, \ox\ra$ for all $x\in \dom f$, and $\la v, x_0\ra<\la v, \ox\ra$ with $x_0\in \dom f$. Thus, the set $\dom f$ and $\{x_0\}$ can be properly separated, so $x_0\notin \mbox{\rm ri}(\dom f)$, which is a contradiction. Therefore, dividing both sides of \eqref{e} by $\gamma$ and use this inequality with $x\in \dom f$ and $\lambda=f(x)$, we see that $v/\gamma\in \partial f(\ox)$.

Now, suppose that $f$ is continuous at $x_0$. Then $x_0\in \mbox{\rm int}(\dom f)$, so $x_0\in \mbox{\rm ri}(\dom f)=\mbox{\rm int}(\dom f)$. Therefore, $\partial f(\ox)\neq\emptyset$. \\[1ex]
(c) Suppose that $f$ is continuous at $\ox$ and that $f(\ox)<\bar{\lambda}$.  By Corollary \ref{intepi}, we see that $(\ox, \bar{\lambda})\in \mbox{\rm int}(\epi f)$ and thus $N((\ox, \bar{\lambda}); \epi f)=\{(0, 0)\}$.

The proofs of the last two assertions are quite obvious, so we left them for the reader.$\h$

The following corollary is a direct consequence of Theorem \ref{CSR}; see also \cite[Theorem~4.29.]{LM2019}.

\begin{Corollary} Let $f_i\colon \R^n\to \oR$ for $i=1, \ldots, m$ be proper nearly convex functions. Suppose that
\begin{equation*}
\bigcap_{i=1}^m \mbox{\rm ri}(\dom f_i)\neq\emptyset.
\end{equation*}
Then $f_1+\cdots+f_m$ is nearly convex and we have the equality
\begin{equation*}
\partial (f_1+\cdots +f_m)(\ox)=\partial f(\ox)+\cdots +\partial f_m(\ox)\ \, \mbox{\rm for all }\, \ox\in \bigcap_{i=1}^m \dom f_i.
\end{equation*}
\end{Corollary}
{\bf Proof.} We only need to prove the result for the case where $m=2$. It suffices to apply Theorem \ref{CSR} to the epigraphical mappings $F_i=E_{f_i}$ for $i=1, 2$ from \eqref{Emapping} with the use of Proposition \ref{CS}.  $\h$

Now we consider the  composition of two mappings $F\colon \R^n\tto \R^p$ and $G\colon \R^p\tto \R^q$. Given $\oz\in(G\circ F)(\ox)$, we consider the set
\begin{equation*}
	M(\ox,\oz)=F(\ox)\cap G^{-1}(\oz).
\end{equation*}
The following theorem provides the coderivative chain rule for nearly convex set-valued mappings.

\begin{Theorem}\label{scr} Let $F\colon \R^n\tto \R^p$ and $G\colon \R^p\tto \R^q$ be  nearly convex set-valued mappings. Suppose that the qualification condition \eqref{QCC} is satisfied.  Then for any $(\ox,\oz)\in\mbox{\rm gph}(G\circ F)$ and $ w\in \R^q$ we have the coderivative chain rule
	\begin{equation*}
		D^*(G\circ F)(\ox,\oz)(w)=\big(D^*F(\ox,\oy)\circ D^*G(\oy,\oz)\big)(w)
	\end{equation*}
	whenever $\oy\in M(\ox,\oz)$.
\end{Theorem}
{\bf Proof.}  Picking $u\in D^*(G\circ F)(\ox,\oz)( w)$ and $\oy\in M(\ox,\oz)$ gives us the inclusion $$(u,- w)\in N((\ox,\oz);\mbox{\rm gph}(G\circ F)).$$  Consider the sets $\Omega_i$ for $i=1, 2$ given in \eqref{O12c}.
We can directly deduce from  the definition of the normal cone that
\begin{equation*}
	(u,0,-w)\in N((\ox,\oy,\oz);\Omega_1\cap\Omega_2).
\end{equation*}
By the proof of Theorem \ref{thm:chain_rule_setvalued} we have $\ri\Omega_1\cap \ri\Omega_2\neq\emptyset$.  Applying Theorem~\ref{NCI}  gives us
\begin{equation*}
	(u,0,-w)\in N((\ox,\oy,\oz);\Omega_1\cap\Omega_2)=N((\ox,\oy,\oz);\Omega_1)+N((\ox,\oy,\oz);\Omega_2).
\end{equation*}
Then the rest of the proof follows that of \cite[Theorem~3.38]{bmn}. $\h$

The next result is a direct consequence of Theorem \ref{scr}.

\begin{Corollary} Let $B\colon \R^n\to \R^p$ be the affine mapping defined by
\begin{equation*}
B(x)=Ax+b\; \mbox{\rm for }x\in \R^n,
\end{equation*}
where $A\in \R^{p\times n}$ and $b\in \R^p$, and let $g\colon \R^p\to \oR$ be a proper nearly convex function. Suppose that
\begin{equation*}
B(\R^n)\cap \mbox{\rm ri}(\dom g)\neq\emptyset.
\end{equation*}
Then  we have the equality
\begin{equation*}
\partial (g\circ B)(\ox)=A^T\partial g(B(\ox))\; \mbox{\rm for every }\ox\in \mbox{\rm dom}(g\circ B).
\end{equation*}
\end{Corollary}
{\bf Proof.} This is a direct consequence of Theorem \ref{scr} with $F(x)=\{B(x)\}$ for $x\in \R^n$, and $G(y)= [g(y), \infty)$ for $y\in \R^p$. $\h$

The next result is another direct consequence of Theorem \ref{scr}.
\begin{Corollary}
	Let $F\colon \R^n\tto\R^p$ be a nearly convex set-valued mapping and let $\Theta\subset\R^p$ be a nearly convex set. Suppose that the qualification condition \eqref{QCIN} is satisfied.  Then we have the equality
\begin{equation*}
N(\ox; F^{-1}(\Theta))=D^*F(\ox, \oy)(N(\oy; \Theta)
\end{equation*}
whenever $\ox\in F^{-1}(\Theta)$ and $\oy\in F(\ox)\cap \Theta$.
\end{Corollary}

Next, we discuss coderivatives of the intersection mapping.

\begin{Theorem}\label{SVMI}  Let $F_i\colon \R^n\tto \R^p$ for $i=1, \ldots, m$ be nearly convex set-valued mappings, and let $F=\bigcap_{i=1}^mF_i$. Assume that the qualification condition \eqref{QCM} is satisfied. Then for any $(\ox, \oy)\in \gph F$ we have
\begin{equation}\label{CIM}
D^*F(\ox, \oy)(y^*)=\bigcup\big\{D^*F_1(\ox, \oy)(y^*_1)+\cdots +D^*F_m(\ox, \oy)(y^*_m)\; \big|\; y^*=y^*_1+\cdots+y^*_m\big\}.
\end{equation}
\end{Theorem}
{\bf Proof.} By Theorem~\ref{NCI}, the set-valued mapping $\bigcap_{i=1}^mF_i$ is nearly convex. Take any $x^*\in D^*F(\ox, \oy)(y^*)$, where $(\ox, \oy)\in \gph F$. Then $(x^*, -y^*)\in N((\ox, \oy); \gph F)$. Since $\gph F=\bigcap_{i=1}^m \gph F_i$, by Theorem~\ref{NCI} we have
\begin{equation*}
(x^*, -y^*)\in N((\ox, \oy); \gph F)=N((\ox, \oy); \gph F_1)+\cdots +N((\ox, \oy); \gph F_m).
\end{equation*}
Thus, there exist $x^*_1, \ldots, x^*_m\in \R^n$ and $y^*_1, \ldots, y^*_m\in \R^p$ such that
\begin{equation*}
x^*=\sum_{i=1}^m x^*_i, \; y^*=\sum_{i=1}^m y^*_i, \; (x^*_i, -y^*_i)\in N((\ox, \oy); \gph F_i).
\end{equation*}
By the definition, $x^*_i\in D^*F_i(\ox, \oy)(y^*_i)$ and hence
\begin{equation*}
x^*\in D^*F_1(\ox, \oy)(y^*_1)+\cdots +D^*F_m(\ox, \oy)(y^*_m).
\end{equation*}
This justifies the inclusion $\subset$ in \eqref{CIM}. The reverse inclusion follows directly from the definition. $\h$

\begin{Corollary} Let $f_i\colon \R^n\to \oR$ be proper nearly convex functions. Suppose that all functions $f_i$  are continuous at $\ox\in \R^n$. Then the maximum function $f$ defined in \eqref{MF} is nearly convex and we have
\begin{equation}\label{sub max}
\partial f(\ox)=\co \big[\bigcup_{i\in I(\ox)}\partial f_i(\ox)\big],
\end{equation}
where $I(\ox)=\{i=1, \ldots, m\; |\; f_i(\ox)=f(\ox)\}$.
\end{Corollary}
{\bf Proof.} The near convexity of $f$ follows from Corollary \ref{NCMf} since $\ox\in \mbox{\rm int}(\dom f_i)=\mbox{\rm ri}(\dom f_i)$ for all $i=1, \ldots, m$ under the continuity of $f_i$ at $\ox$. Consider the epigraphical mappings $F_i=E_{f_i}$ defined in the proof of Corollary \ref{NCMf}. Then
\begin{equation*}
E_f=\bigcap_{i=1}^m F_i.
\end{equation*}
Fix any $v\in \partial f(\ox)$ and get by Proposition \ref{CS} that $v\in D^*E_f(\ox, \bar{\lambda})(1)$, where $\bar{\lambda}=f(\ox)$. By Theorem \ref{SVMI}, there exist $v_i\in \R^n$ and $\lambda_i\in \R$ for $i=1, \ldots, m$ such that
\begin{equation*}
v=\sum_{i=1}^m v_i, \; \sum_{i=1}^m\lambda_i=1, \; (v_i, -\lambda_i)\in N((\ox, \bar{\lambda}); \epi f_i).
\end{equation*}
By Proposition \ref{sub pro}(c), the continuity of $f_i$ at $\ox$ ensures that $ N((\ox, \bar{\lambda}); \epi f_i)=\{(0, 0)\}$ if $i\notin I(\ox)$, i.e., $\bar{\lambda}>f_i(\ox)$. Thus,
\begin{equation*}
v=\sum_{i\in I(\ox)}v_i,\;  \sum_{i\in I(\ox)} \lambda_i=1,
\end{equation*}
where $(v_i, -\lambda_i)\in N((\ox, f_i(\ox)); \epi f_i)$ for $i\in I(\ox)$. Then $\lambda_i\geq 0$ whenever $i\in I(\ox)$ by Proposition \ref{sub pro} (a). In addition,  if $\lambda_i=0$, then by Proposition \ref{sub pro}(d) we have $v_i\in N(\ox; \dom f_i)=\{0\}$ since $\ox\in \mbox{\rm int}(\dom f_i)$. By Proposition \ref{sub pro}(b), we see that $\partial f_i(\ox)\neq\emptyset$ for every $i$ under the continuity of $f_i$. We can also see from Proposition \ref{sub pro}(e) that if $\lambda_i>0$, then $v_i\in \lambda_i\partial f_i(\ox)$ for $i\in I(\ox)$. Then we have $v_i\in \lambda_i\partial f_i(\ox)$ for all $i\in I(\ox)$. Thus,
\begin{equation*}
v\in \sum_{i\in I(\ox)}\lambda_i\partial f_i(\ox), \; \lambda_i\geq 0, \; \sum_{i\in I(\ox)}\lambda_i=1.
\end{equation*}
This justifies the inclusion $\subset$ in \eqref{sub max}. The reverse inclusion follows directly from the definition. $\h$\\[1ex]
{\bf Acknowledgment.} Nguyen Mau Nam would like to thank the Vietnam Institute of Mathematics-VAST (through the IM-Simons program) and the Vietnam Institute for Advanced Study in Mathematics for hospitality.


\vspace*{-0.2in}

\small
\begin{thebibliography}{99}


\bibitem{nr1} H. H. Bauschke, W. L. Hare, W. M. Moursi,  On the range of the
Douglas-Rachford operator. {\em Math. Oper. Res.} {\bf 41} (2016), no. 3, 884--897.

\bibitem{bmw2013} H. H. Bauschke, S. M. Moffat and X. Wang, Near equality, near convexity, sums of maximally monotone operators, and averages of firmly nonexpansive mappings. {\em Math. Program.} \textbf{139} (2013), Ser. B, 55--70.


\bibitem{nr2}  H. H. Bauschke, W. M.  Moursi,  On the behavior of the Douglas-Rachford
algorithm for minimizing a convex function subject to a linear constraint. {\em SIAM J. Optim.} {\bf 30} (2020), no.
3, 2559--2576.
%
\bibitem{Borwein2000} J. M. Borwein and A. S. Lewis, {\em Convex Analysis and Nonlinear Optimization}, 2nd edition, Springer, New York, 2006.

\bibitem{bkw2008} R. I. Bo{\c{t}}, G. Kassay  and G. Wanka, Duality for almost convex optimization problems via the perturbation approach, {\em J. Global Optim.} \textbf{42} (2008), 385–399.

\bibitem{nr3} D. Hinrichsen, E. Oeljeklaus,  The set of controllable multi-input systems is
generically convex. {\em Math. Control Signals Systems.} {\bf 31} (2019), no. 3, 265--278.

\bibitem{HU2} J. B. Hiriart-Urruty and C. Lemar\'echal, {\em Convex Analysis and Minimization Algorithms, I: Fundamentals}, Springer, Berlin, 1993.


\bibitem{LM2019}    J. Li and G. Mastroeni, Near equality and almost convexity of functions with applications to optimization and error bounds, {\em J. Convex Anal.} {\bf 26} (2019), no. 3, 785--822.


    \bibitem{nr4} W. M. Moursi, The forward-backward algorithm and the normal problem. {\em J.
Optim. Theory Appl.} {\bf  176} (2018), no. 3, 605--624.

    \bibitem{nr5} H. Luo, X. Wang and B. Lukens,
Variational analysis on the signed distance functions, {\em J.
Optim. Theory Appl.} {\bf 180} (2019), 751--774.




\bibitem{Minty1961} G. J. Minty, On the maximal domain of a ”monotone” function, Michigan Math. J. 8,
pp. 135–137, 1961.


\bibitem{bmn} B. S. Mordukhovich and N. M. Nam, {\em An Easy Path to Convex Analysis and Applications: Second Edition}, Springer, Forthcoming.

\bibitem{bmn2022}  B. S. Mordukhovich and N. M. Nam, \textit{Convex Analysis and Beyond, Vol. I: Basic Theory}, Springer, Cham, Switzerland, 2022.

\bibitem{mmw2016} S. M. Moffat, W. M. Moursi and X. Wang,
Nearly convex sets: fine properties and domains or ranges of subdifferentials of convex functions,
{\em Math. Program.} \textbf{160} (2016), Ser. A, 193--223.


\bibitem{r} R. T. Rockafellar, {\em Convex Analysis}, Princeton University Press, Princeton, NJ, 1970.


\bibitem{R1970} R.T. Rockafellar, On the virtual convexity of the domain and range of a nonlinear
maximal monotone operator, Math. Ann. 185, pp. 81–90, 1970.


\bibitem{rw} R. T. Rockafellar and R. J-B.
Wets, {\em Variational Analysis}, Springer,  Berlin,
1998.


\end{thebibliography}
\end{document}